\documentclass[11 pt,Rectoverso]{article}  
\newcommand{\indic}{1\negthickspace\textrm{I}} 

\usepackage[latin1]{inputenc}
\usepackage{graphicx,psfrag,amsfonts,rotating}
\usepackage{amsmath,indentfirst,floatflt}
\usepackage{qsymbols}
\def \HH{{\cal{H}}}
\def \W{{\cal W}}
\def \R{{\bf R}}
\def \B{{\cal B}}
\def \E{{\cal E}}
\def \M{{\cal M}}
\def \S{{\bf S}}
\def \Z{{\bf Z}}
\def \H{{\bf H}}
\def \e{{\sf e}}
\def \b{{\sf b}}
\def \w{{\sf w}}
\def \m{{\sf m}}
\def \x{{\sf x}}
\def \p{{\sf p}}

\def \tb{{\tiny $\blacksquare$~ }}
\def \app#1#2#3#4#5{\begin{array}{rccl} #1:&#2&\longrightarrow&#3\\ &#4&\longmapsto&#5\end{array}}

\def \k#1{{\cal K}_{#1}}

\def \be{\begin{eqnarray*}}
\def \ee{\end{eqnarray*}}
\def \ben{\begin{eqnarray}}
\def \een{\end{eqnarray}}
\def \bq{\begin{equation}}
\def \eq{\end{equation}}

\def \build#1#2#3{\mathrel{\mathop{\kern 0pt#1}\limits_{#2}^{#3}}}
\def \cro#1{\llbracket#1\rrbracket}

\def \ceil#1{\lceil#1\rceil}
\def \floor#1{\lfloor#1\rfloor}
\def\proof{\noindent{\bf Proof. }}

\def \imp{\Rightarrow}

\def \eref#1{(\ref{#1})}
\def \sous#1#2{\mathrel{\mathop{\kern 0pt#1}\limits_{#2}}}
\def \sur#1#2{\mathrel{\mathop{\kern 0pt#1}\limits^{#2}}}

\def \g#1{\big[#1\big]}
\def \dd{\xrightarrow[n]{(d)}}

\def \captionn#1{\begin{center}\begin{minipage}{14cm}\sf\caption{\small #1}\end{minipage}\end{center}}

\def \dis{\displaystyle}
\def \tend{\longrightarrow}

\def\l{\left}
\def\r{\right}

\def \bar{\overline}

\def \u{u}
\def \d{d}
\def \tu{{\mathop{T}^{(\u)}}}
\def \tuu{{\mathop{T}^{(\u\u)}}}
\def \tud{{\mathop{T}^{(\u\d)}}}
\def \tdu{{\mathop{T}^{(\d\u)}}}
\def \tdd{{\mathop{T}^{(\d\d)}}}
\def \td{{\mathop{T}^{(\d)}}}
\def \tu{{\mathop{T}^{(\u)}}}
\def \bij{\leftrightarrow}
\def \vc{\varcurlyvee}
\def \QQ{\widetilde{Q}_n} 
\font\dsrom=dsrom10 scaled 1200

\def \ind{\textrm{\dsrom{1}}}

\newqsymbol{`P}{\mathbb{P}}
\newqsymbol{`R}{\mathbb{R}}
\newqsymbol{`Z}{\mathbb{Z}}
\newqsymbol{`E}{\mathbb{E}}
\newqsymbol{`e}{\varepsilon}
\newqsymbol{`a}{\alpha}
\newqsymbol{`t}{\tau}
\newqsymbol{`O}{\Omega}
\newqsymbol{`M}{\overrightarrow{\cal M}}
\newqsymbol{`m}{{M_n^{\bullet}}}
\newqsymbol{`Q}{{\mathbb{Q}}}
\newqsymbol{`q}{{{\cal Q}_n^{\bullet}}}
\newqsymbol{`W}{{\cal W}^+}
\newqsymbol{`w}{{\cal W}}
\newqsymbol{`d}{{\mathcal D}}
\newqsymbol{`g}{{\mathcal G}}

\newqsymbol{`N}{\mathbb{N}}

\graphicspath{{../DESSINS/}}
\renewcommand{\baselinestretch}{1.2}

\DeclareMathOperator{\Ref}{Ref}
\DeclareMathOperator{\Cont}{Cont}

\begin{document}
\newtheorem{lem}{Lemma}
\newtheorem{defi}{Definition}
\newtheorem{pro}[lem]{Proposition}
\newtheorem{theo}[lem]{Theorem}
\newtheorem{cor}[lem]{Corollary}
\newtheorem{remi}{Remark\rm}{\rm}
\newtheorem{comm}{Comments\rm}{\rm}
\newenvironment{rem}%
{\begin{center}\begin{minipage}{16cm}\begin{remi}}%
{\end{remi}\end{minipage}\end{center}}

\begin{center}
\Large{\bf  Asymptotics of Bernoulli random walks, bridges, excursions and meanders with a given number of peaks }
\large
{\[\begin{array}{lcl}
\textsf{Jean-Maxime Labarbe }&~~& \textsf{Jean-François Marckert} \\
\textrm{Universit\'e de Versailles}&~~ & \textrm{CNRS, LaBRI}\\
\textrm{45 avenue des Etats-Unis}&~~ &\textrm{Université Bordeaux 1}\\
\textrm{78035 Versailles cedex}&~~& \textrm{351 cours de la Libération}\\
\textrm{labarbe@math.uvsq.fr} & ~~ & \textrm{33405 Talence cedex}\\
&~~ &\textrm{marckert@labri.fr}
 \end{array}\]
}

\end{center}

\begin{abstract}A Bernoulli random walk is a random trajectory starting from 0 and having i.i.d. increments, each of them being $+1$ or $-1$, equally likely. The other families cited in the title are Bernoulli random walks under various conditionings. A peak in a trajectory is a local maximum. In this paper, we condition the families of trajectories to have a given number of peaks. We show that, asymptotically, the main effect of setting the number of peaks is to change the order of magnitude of the trajectories. The counting process of the peaks, that encodes the repartition of the peaks in the trajectories, is also studied. It is shown that suitably normalized, it converges to a Brownian bridge which is independent of the limiting trajectory. Applications in terms of plane trees and parallelogram polyominoes are also provided.
\end{abstract}

\section{Introduction}
Let $`N=\{0,1,2,3,\dots\}$ be the set of non-negative integers.
For any $n\in `N$, we denote by $\W_n$ the set of Bernoulli chains with $n$ steps~:
\[\W_n=\{ {\S}=(S(i))_{0\leq i \leq n}~: S(0)=0, S({i+1})=S(i) \pm 1 \textrm{ for any }i\in\cro{0,n-1}\}.\]
The sets of Bernoulli bridges $\B_n$,  Bernoulli excursions $\E_n$, Bernoulli meanders $\M_n$ with $n$ steps are defined by
\be
\B_n&=&\{ \S~: \S\in\W_n, S(n)=0\},\\
\E_n&=&\{\S~:\S\in\W_n, S(n)=0, S(i)\geq 0 \textrm{ for any }i \in\cro{0,n}\},\\
\M_n&=&\{\S~:\S\in\W_n, S(i)\geq 0 \textrm{ for any }i \in\cro{0,n}\}. 
\ee  
The cardinalities of these sets are given by
\begin{equation}\label{card2}
\#\W_n=2^n,~~~ 
\#\B_{2n}=\binom{2n}{n},~~~~ \#\E_{2n}=\frac{1}{n+1}{\binom{2n}{n}},~~~ \#\M_n=\binom{n}{\floor{n/2}}, 
\end{equation} and
for every odd number $n$, $\B_n=\E_n=\varnothing$. 
The two first formulas are obvious, the third can be proved for instance thanks to the cyclical lemma (see also the 66 examples of the appearance of the Catalan numbers $\#\E_{2n}$ in combinatorics in Stanley \cite[ex. 6.19 p.219]{STAN}), and the last one, may be proved iteratively or thanks to a bijection with Bernoulli bridges (see Section \ref{bml}).

Let $n\in\mathbb{N}$. For every $\S\in \W_n$, the set of peaks of $\S$, denoted by $\S_{\wedge}$, is defined by
\[\S_{\wedge}=\{x~: x\in\cro{1,n-1},~ S({x-1})=S({x+1})=S(x)-1\}.\]
The set $(-\S)_{\wedge}$ is called the set of valleys of $\S$~: it is easy to check that for any $\S$,  $\#\S_{\wedge}-\#(-\S)_{\wedge}$ belongs to $\{+1,0,-1\}$. The value of this difference depends only on the signs of the first and last steps of $\S$. In this paper, we focus only on the number of peaks and we denote by $\W_n^{(k)}$ (resp. $\B_n^{(k)}$, $\E_n^{(k)}$  and $\M_n^{(k)}$)  the subset of $\W_n$, (resp. $\B_n$, $\E_n$, $\M_n$) of trajectories having exactly $k$ peaks (for any $k> \floor{n/2}$ these sets are empty). 
We have
\begin{pro}\label{card} For any $k\geq 0$ and any $n\geq 0$, 
\[
\#\W_n^{(k)} =  \binom{n+1}{n-2k},~
\#\B_{2n}^{(k)} = \binom{n}{k}^2,~
\#\E_{2n}^{(k)} = \frac{1}{n}\binom{n}{k}\binom{n}{k-1},~
\#\M_{n}^{(k)} = \binom {\floor {n/2}} k\binom {\ceil {n/2}} k,  \]
where, by convention
\[\binom mp=\left\{
 \begin{array}{l}
\dis\frac{m!}{p!(m-p)!} \textrm{ if }p \textrm{ and }m \textrm{ are non negative integers and } p\in\cro{0,m},\\
0 \textrm{ in any other cases.}
 \end{array}
\right.\]
\end{pro}
\begin{figure}[htbp]
\centerline{\includegraphics[height=4cm]{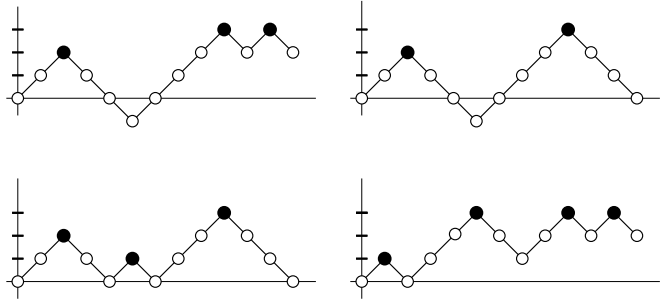}}
\caption{\label{exe}Trajectories from ${\cal W}_{12}^{(3)}$, ${\cal B}_{12}^{(2)}$, ${\cal E}_{12}^{(3)}$, and ${\cal M}_{12}^{(4)}$. Black dots correspond to peaks.}
\end{figure}
The formula giving $\#\E_{2n}^{(k)}$ is due to Narayana \cite{NAR1} computed in relation with pairs of $k$-compositions of $n$ satisfying some constraints (see also Stanley \cite[ex. 6.36 p.237]{STAN}). \par

Let $`P_n^{\w}$, $`P_n^{\b}$, $`P_n^{\e}$  and $`P_n^{\m}$ be the uniform law on $\W_n$, $\B_n$,  $\E_n$, and $\M_n$ and $`P_n^{\w,(k)}$, $`P_n^{\b,(k)}$, $`P_n^{\e,(k)}$  and $`P_n^{\m,(k)}$  be the uniform law on $\W_n^{(k)}$, $\B_n^{(k)}$, $\E_n^{(k)}$ and $\M_n^{(k)}$. For $\x\in\{\w,\b,\e,\m\}$, a random variable under $`P_n^{\x,(k)}$ is then simply a random variable under $`P_n^{\x}$  conditioned to have $k$ peaks.
We are interested in the asymptotic behavior of random chains under the distributions $`P_n^{\x,(k)}$, when $n$ and $k=\k{n}$  go to infinity.
\par
Let $C[0,1]$ be the set of continuous functions defined on $[0,1]$ with real values. For any $\S\in \W_n$, denote by $u_n$ the function in $C[0,1]$ obtained from $\S$ by the following rescaling: 
\begin{equation}\label{uuuu}
u_n(t)=\frac{1}{\sqrt{n}}\big(S({\floor{nt}})+\{nt\}(S({\ceil{nt}})-S({\floor{nt}}))\, \big)\textrm{ for any }t\in[0,1].
\end{equation}

We call Brownian bridge $\b$, Brownian excursion $\e$ and Brownian meander $\m$ the (normalized) processes characterized as follows~: let $\w$ be a 1-dimensional standard Brownian motion. Let $d=\inf\{t~: t\geq 1, \w_t=0\}$ and $g=\sup\{t~: t\leq 1, \w_t=0\}$. Almost surely, we have $g-d>0$, $g\in(0,1)$.
The processes $\b,\e$ and $\m$ have the following representations~:
\[(\b_t)_{t\in[0,1]}\sur{=}{(d)} \Big(\frac{\w_{gt}}{\sqrt{g}}\Big)_{t\in[0,1]},~~ (\e_t)_{t\in[0,1]}\sur{=}{(d)} \Big(\frac{|\w_{d+(g-d)t}|}{\sqrt{d-g}}\Big)_{t\in[0,1]},~~  (\m_t)_{t\in[0,1]}\sur{=}{(d)}  \Big(\frac{|\w_{(1-g)t}|}{\sqrt{1-g}}\Big)_{t\in[0,1]}.\]
As a consequence of the Donsker \cite{BIL} theorem (for $\x=\w$), Kaigh \cite{KWD} (for $\x=\b$), Iglehart \cite{IGL} (for $\x=\m$) and Belkin \cite{BEL} (for $\x=\e$), 
\begin{theo}\label{dgh} For any $\x\in\{\w,\b,\e,\m\}$,
under $`P_n^{\,\x}$, $u_n\dd\x$ in $C[0,1]$ endowed with the topology of the uniform convergence.
\end{theo}
In the case $\x\in\{\b,\e\}$, even if non specified, it is understood that $n\to +\infty$ in $2\mathbb{N}$. \par
In fact, Theorem \ref{dgh} can be proved directly, thanks to the elementary enumeration of paths passing via some prescribed positions in the model of Bernoulli paths. The method used to show the tightnesses in our Theorem \ref{pri} may be used to prove the tightness in Theorem \ref{dgh}; thanks to some probability tricks, this reduces to show the tightness under $`P_n^{\,\w}$, which is simple.  
\medskip

The finite dimensional distributions of $\w,\e,\b$ and $\m$ are recalled in Section \ref{fdd}. Numerous relations exist between these processes, and their trajectories, and a lot of parameters have been computed. We refer to Bertoin \& Pitman \cite{BP}, Biane \& Yor \cite{BY}, Pitman \cite{PIT14} to have an overview of the subject. These convergences have also provided some discrete approaches to the computation of values attached to these Brownian processes, and the literature about that is considerable, see e.g. Cs{á}ki \& Y.  Hu \cite{CH}, and reference therein.\par
We introduce the counting process of the number of peaks~: for any $\S\in{\W_n}$, denote by $\Lambda(\S)=(\Lambda_l(\S))_{l\in \cro{0,n}}$ the process~:
\begin{equation}
\Lambda_l(\S)=\#\S_{\wedge}\cap\cro{0,l} \textrm{ for any }l\in \cro{0,n}.
\end{equation}
For $\S\in\W_n$, $\Lambda_n(\S)=\#\S_{\wedge}$ is simply the total number of peaks in $\S$. 
We have
\begin{pro} \label{CVD}For any $\x\in\{\w,\b,\e,\m\}$, under $`P_n^\x$,
\[\frac{\Lambda_n-n/4}{\sqrt{n}}\xrightarrow[~n~]{(d)}\mathcal{N}(0,1/16),\]
where $\mathcal{N}(0,1/16)$ denotes the centered Gaussian distribution with variance $1/16$.
\end{pro}

We will now describe the main result of this paper. Its aim is to describe the influence of the number of peaks on the shape of the families of trajectories introduced above. We will then condition the different families by $\#\S_{\wedge}=\k{n}$ for a general sequence $(\k{n})$ satisfying the following constraints~:
\[(H)=\left(
\hbox{For any } n,~~ \k{n}\in \mathbb{N},
~~\dis\lim_n \k{n}= +\infty, ~~\dis\lim_n n/2-\k{n}= +\infty
\right).\]
Notice that for every $\S\in\W_n$, $\#\S_{\wedge}\in\cro{0,\floor{n/2}}$ and then $(H)$ is as large as possible to avoid that the sequences $\k{n}$ and $n/2-\k{n}$ have a finite accumulation point.

We set $p_n:={2\k{n}}/{n}$ and 
\begin{equation}\label{beb}
\beta_n:=\sqrt{n(1-p_n)/p_n},\hbox{ and } \gamma_n:= \sqrt{np_n(1-p_n)}= p_n\beta_n.
\end{equation}
Each peak can be viewed to be made by two consecutive steps; hence, if you pick at random one step of a trajectory under $`P_n^{\x,(\k{n})}$,  the probability that this step belongs to a peak is $p_n$.

We consider $\S$ and $\Lambda(\S)$ as two continuous processes on $[0,n]$, the values between integer points being defined by linear interpolation. The normalized versions of $\S$ and $\Lambda(\S)$ are respectively denoted by $s_n$ and $\lambda_n$~:
\begin{equation}\label{sl}
s_n(t):=\frac{S(nt)}{\beta_n}~~~~\textrm{ and }~~~ \lambda_n(t):={2}\,\,\frac{\Lambda_{nt}-t\,\k{n}}{\gamma_n}.\end{equation}
\begin{theo}\label{pri} If $(H)$ is satisfied, for any $\x\in\{\w, \b, \e,\m\}$, under $`P_n^{\,\x,(\k{n})}$, 
\begin{equation}
(s_n,\lambda_n)\dd (\x,\widehat{\b})
\end{equation}
where $\widehat{\b}$ is a Brownian bridge independent of $\x$ and where the weak convergence holds in $C([0,1])^2$ endowed with the topology of uniform convergence.
\end{theo}
Hence, under $`P_n^{\,\x,(\k{n})}$, up to the scaling constant, the process $s_n$ behaves as under $`P_n^{\,\x}$. The normalizing factor $\beta_n$, that will be explained later in the paper, indicates the order of magnitude of the process $\S$ under $`P_n^{\,\x,(\k{n})}$ ($\beta_n$ is a decreasing function of $\k{n}$). The normalizing constant $\gamma_n$ is smaller than $\sqrt{n/4}$ whatever is $p_n$; $\gamma_n$ gives the asymptotic order of the ``linearity defect'' of $t\mapsto \Lambda_{nt}$.
The fact that $(\lambda_n)$ converges to a Brownian bridge independent of the limit trajectory is quite puzzling. For example, under  $`P_n^{\,\e}$, one would expect that only few peaks appear in a neighborhood of 0, this resulting in a negative bias in $\lambda_n$ near 0. This must be true, but this bias is not important enough to change the asymptotic behavior of $\lambda_n$. \par

A second direct corollary of Theorem \ref{pri} is stated below:
\begin{cor}\label{zeze} For any $\x\in\{\w, \b, \e,\m\}$, under $`P_n^{\,\x}$, we have
\begin{equation}\label{ww}
\left(\frac{S({nt})}{\sqrt{n}}~,~4\,\,\frac{\Lambda_{nt}-tn/4}{\sqrt{n}}\right)_{t\in[0,1]}\dd (\x_t,\widehat{\w}_t)_{t\in[0,1]}
\end{equation}
where $\widehat{\w}$ is a Brownian motion independent of $\,\x$ and where the weak convergence holds in $C([0,1])^2$ endowed with the topology of uniform convergence.
\end{cor} 
Theorem \ref{dgh} is of course a consequence of Corollary \ref{zeze}.\medskip

\proof For any $\S\in\W_n$, set $q_n(\S)=2\Lambda_n(\S)/n$, $\tilde{\beta}_n(\S)= \sqrt{n(1-q_n(\S))/q_n(\S)},$ and $\tilde{\gamma_n}(\S):= \sqrt{nq_n(\S)(1-q_n(\S))}$. We have
\begin{equation}\label{fi}
\left(\frac{S(nt)}{\sqrt{n}},~4\frac{\Lambda_{nt}-tn/4}{\sqrt{n}}\right)=\left(\frac{S(nt)}{\tilde{\beta}_n}\frac{\tilde{\beta}_n}{\sqrt{n}},~2\frac{\tilde{\gamma_n}}{\sqrt{n}}\l(2\,\frac{\Lambda_{nt}-t\Lambda_n}{{\tilde{\gamma_n}}}\r)+4t\,\frac{\Lambda_n-n/4}{\sqrt{n}}\right).\end{equation}
By Proposition \ref{CVD} and Theorem \ref{pri}, under $`P_n^{\x}$, the five-tuple
\[\left(4\,\frac{\Lambda_n-n/4}{\sqrt{n}},\, \l(\frac{S(nt)}{\tilde{\beta}_n}\r)_{t\in[0,1]},2\l(\frac{\Lambda_{nt}-t\Lambda_n}{\tilde{\gamma_n}}\r)_{t\in[0,1]},2\frac{\tilde{\gamma_n}}{\sqrt{n}}, \frac{\tilde{\beta_n}}{\sqrt{n}}\right)\]
converges in distribution to
\[(N,(s_t)_{t\in[0,1]}, (\lambda_t)_{t\in[0,1]},A,B)\]
where $N$ is a centered Gaussian random variable with variance 1, and where, conditionally on $N$, $(s,\lambda)\sur{=}{(\d)}(\x,\widehat{b})$ where $\x$ and $\widehat{b}$ are independent, $\widehat{b}$ is a Brownian bridge, and $A$ and $B$  are two random variables equal to 1 a.s.. By \eref{fi}, $\left(\frac{S(nt)}{\sqrt{n}},4\,\frac{\Lambda_{nt}-tn/4}{\sqrt{n}}\right)$ converges to $(\x,(\widehat{b}_t+tN)_{t\in[0,1]})$ where $N$ is independent of $\x$ and $\widehat{b}$, and then the result follows, since $(\widehat{b}_t+tN)_{t\in[0,1]}$ is a standard Brownian motion.~$\Box$
 
\subsubsection*{Consequences in terms of plane trees}
\begin{figure}[htbp]
\centerline{\includegraphics[height=3cm]{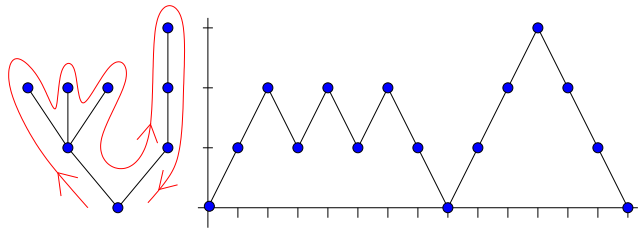}}
\caption{\label{contour}A plane tree and its contour process}
\end{figure}

Consider the set $T_n$ of \it plane trees \rm (rooted ordered trees) with $n$ edges (we refer to \cite{ALD,MM2} for more information on these objects). 
There exists a well known bijection between $T_n$ and $\E_{2n}$ which may be informally described as follows. Consider a plane tree $\tau\in T_n$ (see Figure \ref{contour}), and a fly walking around the tree $\tau$ clockwise, starting from the root, at the speed 1 edge per unit of time. Let $V(t)$ be the distance from the root to the fly at time $t$. The process $V(t)$ is called in the literature, the contour process or the Harris' walk associated with $\tau$. 
The contour process is the very important tool for the study of plane trees and their asymptotics and we refer to Aldous \cite{ALD}, Pitman \cite[Section 6]{PITT}, Duquesne \& Le Gall \cite{DG}, Marckert \& Mokkadem \cite{MM2} for considerations on the asymptotics of normalized trees. It is straightforward that the set of trees encoded by $\E_{2n}^{(k)}$ is  the subset of $T_n$ of trees having exactly $k$ leaves (say $T_n^{(k)}$), a leaf being a node without any child. A corollary of Theorem \ref{pri}, is that random plane tree with $n$ edges and $\k{n}$ leaves, converges, normalized by $\beta_n/2$, to the continuum random tree introduced by Aldous \cite{ALD}, which is encoded by $2\e$. The variable $\Lambda_{nt}$ gives the number of leaves visited at time $nt$. By Theorem \ref{pri}, $\sup_{t\in[0,1]}|(\Lambda_{nt}-t\k{n})/n|\sur{\tend}{proba.}0$. This translates the fact that the leaves are asymptotically uniformly distributed on a random tree chosen equally likely in $T_n^{(\k{n})}$.

\subsubsection*{Consequences in terms of parallelogram polyominoes}

\begin{figure}[htbp]
\centerline{\psfrag{c}{$c$}\psfrag{d}{$d$}\psfrag{t}{$t$}\psfrag{l}{$l$}\psfrag{l-(d-c)}{$l-(d-c)$}
\includegraphics[height=3cm]{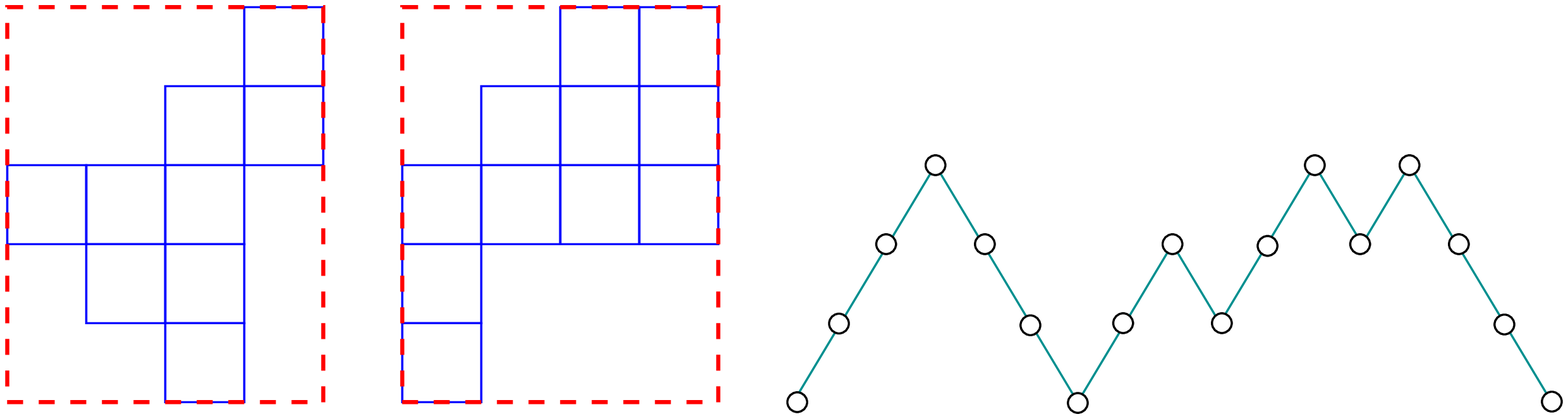}}
\captionn{\label{ppo}The first convex polyomino is not parallelogram, the second is. Their areas are 9 and 11, their perimeters  equal  that of their minimal rectangles, here 18. For both polyominoes $H({\cal P})=4$ and $V({\cal P})=5$. The last picture represents the Bernoulli excursion associated by $\rho$ with the parallelogram polyomino.}
\end{figure}

We refer to Delest \& Viennot \cite{DV} for more information on parallelogram polyominoes. Unit squares having their vertices at integer points in the Cartesian plane are called cells. 
A  polyomino is a  finite union of cells with  connected interior. The number of cells is the area and the length of the border is called the perimeter (see Figure \ref{ppo}). A polyomino ${\cal P}$ is said to be convex if the intersection of ${\cal P}$ with any horizontal or vertical line is a convex segment. For any convex polyomino ${\cal P}$ there exists a minimal rectangle $R({\cal P})$ (that can be seen as a convex polyomino) containing ${\cal P}$. Then ${\cal P}$ touches the border of $R({\cal P})$ along four connected segments. 
A convex polyomino ${\cal P}$ is said to be a parallelogram polyomino if the south-west point and the north-east point of $R({\cal P})$ belongs to ${\cal P}$ (see Figure \ref{ppo}).
Let denote by $H({\cal P})$ and $V({\cal P})$ the horizontal and vertical length of the border of $R({\cal P})$, and let ${\sf Pol}_{n}$ be the set of parallelogram polyominoes with perimeter $n$. 
\begin{pro}(Delest \textrm{\&} Viennot \cite[Section 4]{DV}) For any integer $N \geq 1$, there is a bijection $\rho$ between $\E_{2N}$ and ${\sf Pol}_{2N+2}$, such that if ${\cal P}=\rho(\S)$, the area of ${\cal P}$ is equal to the sum of the heights of peaks of $\S$, moreover $\#\S_{\wedge}=H({\cal P})$, and $V({\cal P})=2N+2-2\#\S_{\wedge}$ (where $2N-2\#\S_{\wedge}$ is the number of steps of $\S$ that do not belong to a peak).
\end{pro}
By symmetry with respect to the first diagonal, the random variables $V({\cal P})$ and $H({\cal P})$ have the same distribution when ${\cal P}$ is taken equally likely in ${\sf Pol}_{2N+2}$. Hence, the proposition says that under $`P_n^{\e}$,  $2N+2-2\#\S_{\wedge}$ and $\#\S_{\wedge}$ have the same distribution. 

We describe in a few words Delest \textrm{\&} Viennot's  bijection: the successive lengths of the columns of the polyomino ${\cal P}$ give the successive heights of the peaks of $\S$. The difference between the heights of the floor-cells of the $i$th and $i+1$th columns of ${\cal P}$ \bf plus one \rm gives the number of down steps between the $i$th and $i+1$th peaks of $\S$. 

For $i\in\{1,\dots,H({\cal P})\}$, let $v_i({\cal P})$ be the number of cells in the $i$th column of ${\cal P}$.  The values $(v_i({\cal P})_{i\in\cro{1,H({\cal P})}}$ coincide with the ordered sequence $(S_i)_{i\in \S_{\wedge}}$. 
Let $`P_{Pol(2n+2)}^{(\k{n})}$ be the uniform law on the set of parallelogram polyominos with perimeter $2n+2$ and width $\k{n}$ (that is such that $H({\cal P})=\k{n})$. Assume that $v$ is interpolated between integer points, and $v(0)=0$.
We have
\begin{pro} If $(\k{n})$ satisfies $(H)$, under $`P_{Pol(2n+2)}^{(\k{n})}$
\[\l(\frac{v(\k{n}t)}{\beta_n}\r)_{t\in[0,1]}\sur{\tend}{(d)} (\e_t)_{t\in[0,1]}\]in $C[0,1]$ endowed with the topology of uniform convergence.
\end{pro}
\proof  Let $(V_i)_{i\in\{1,\dots,\k{n}\}}$ be the successive height of the peaks in $\S$. Assume also that $V(0)=0$ and that $V$ is interpolated between integer points.
By Delest \& Viennot's bijection,  $\beta_n^{-1}v(\k{n}.)$ under $`P_{Pol(2n+2)}^{(\k{n})}$ has the same distribution as 
$\beta_n^{-1}V(\k{n}.)$ under $`P_{n}^{\e,(\k{n})}$. 
Since $(\beta_n^{-1}S({nt}))_{t\in[0,1]}\sur{\tend}{(d)} (\e_t)_{t\in[0,1]}$, to conclude, it suffices to show that 
\begin{equation}\label{zer}
\sup_{t\in[0,1]} \l|\frac{V(\k{n}t)-S({nt})}{\beta_n}\r| \sur{\tend}{proba.}0.
\end{equation}
Let $J(i)$ be  (abscissa of) the $i$th peak in $\S$. 
We have, for any $t \in \{0,1/\k{n},\dots,\k{n}/\k{n}\}$, 
\begin{equation}\label{zer2}V(\k{n}t)-S({nt})=S(J(\k{n}t))-S({nt}).\end{equation}
As one can see using the convergence of $\lambda_n$ to $\widehat{b}$, 
\begin{equation}\label{zer3}
\sup_{t}\l|\frac{J(\k{n}t)-nt}{n}\r|\sur{\tend}{proba.}0.
\end{equation} Indeed, $\sup_t |J(\k{n}t)-nt|/n \leq \sup_t|\Lambda_{nt}-t\k{n}|/n\sur{\tend}{proba.}0$. Since $(s_n)$ converges in $C[0,1]$ under $`P_n^{\e,(\k{n})}$, by a simple argument about its modulus of continuity, using \eref{zer2} and \eref{zer3}, formula \eref{zer} holds true.
$~\Box$\medskip

We would like to point out the work of de Sainte-Catherine \& Viennot \cite{DSV}, who exhibit a quite unexpected link between the cardinalities of excursions having their peaks in a given subset of $\mathbb{N}$ and the famous Tchebichev polynomials.

\section{Combinatorial facts~: decomposition of trajectories}

The decomposition of the trajectories is determinant in our approach, since we will prove directly the convergence of finite dimensional distributions under $`P_n^{\,\x}$. An important difference with the case where the peaks are not considered is that under $`P_n^{\,\x, (k)}$, $\S$ is not a Markov chain. 
Indeed, the law of $(S(j))_{l\leq j\leq n}$ depends on $S(l)$, on the number of peaks in $\cro{0,l}$, and also on the step $(S({l-1}),S(l))$.  The computation of the distributions of vector $(S({t_1}),\dots,S({t_k}))$ in these various cases admits the complication that $\S$ may own a peak in some of the positions $t_1, \dots, t_k$. In order to handle the contribution of these peaks, we have to specify what the type $\u$ or $\d$ ($+1$ or $-1$) of the first and last steps of the studied parts $(S({t_i}),\dots,S({t_{i+1}}))$ is.\par
We set the following notation~: $\Delta S_l=S({l})-S(l-1)$, and write for convenience, $\Delta S_l=\u$ when $\Delta S_l=1$, and $\Delta S_l=\d$ when $\Delta S_l=-1$. In this paper, we deal only with discrete trajectories $\S$ such that $\Delta S_k \in\{+1,-1\}$ for any $k$. We will not recall this condition.~\par
For $a$ and $b$ in $\{\d,\u\}$, and $l,x,y,j$ in $\mathbb{Z}$, set
\be
T_{ab}^{j}(l,x,y)&=&\{\S~:\S=(S(i))_{0\leq i\leq l},~ \#\S_{\wedge}=j,~ \Delta S_1=a,~ \Delta S_l=b,~ S(0)=x,~S(l)=y \}\\
T_{ab}^{j,\geq}(l,x,y)&=&\{\S~: \S \in T_{ab}^{j}(l,x,y),~ S(i)\geq 0 \textrm{ for any }i\in\cro{0,l}\}.
\ee
For any $l,j_1,j_2,x,y\in\mathbb{Z}$, set
\be
\g{l,j_1,j_2,x,y}=\binom{\dis\frac{l+y-x}{2}-1}{j_1}\binom{\dis\frac{l-y+x}{2}-1}{j_2}.\ee
We have
\begin{pro}\label{denu} For any $y\in\mathbb{Z}$, $l\geq 0$, $j\geq 0$, 
\begin{equation}\label{resum1}
\#T_{ab}^{j}(l,0,\,y)=\g{l,\,j-\ind_{b=\d},\,j-\ind_{a=\u},\,0,\,y}.
\end{equation}

 For any $x\geq 0$, $y\geq 1$, $l\geq 0$, $j\geq 0$,  
\be
\#T_{\u\u}^{j,\geq}(l,\,x,\,y)&=&\g{l,\,j,\,j-1,\,x,\,y}-\g{l,\,j-1,\,j,\,-x,\,y},\\
\#T_{\u\d}^{j,\,\geq}(l,\,x,\,y)&=&\g{l,\,j-1,\,j-1,\,x,\,y}-\g{l,\,j-2,\,j,\,-x,\,y},
\ee
For any $x\geq 1$, $y\geq 1$, $l\geq 0$, $j\geq 0$,
\be
\#T_{\d\u}^{j,\,\geq}(l,\,x,\,y)&=&\g{l,\,j,\,j,\,x,\,y}-\g{l,\,j,\,j,\,-x,\,y},\\
\#T_{\d\d}^{j,\,\geq}(l,\,x,\,y)&=&\g{l,\,j-1,\,j,\,x,\,y}-\g{l,\,j-1,\,j,\,-x,\,y}
\ee
(notice that $\#T_{\d b}^{j,\geq}(l,0,y)=0$ and $\#T_{a \u}^{j,\geq}(l,x,0)=0$). In other words
\begin{equation}\label{resum2}
\#T_{ab}^{j,\,\geq}(l,x,\,y)=\g{l,\,j-\ind_{b=\d},\,j-\ind_{a=\u},\,x,\,y}-\g{l,\,j-\ind_{a=\u}-\ind_{b=\d},\,j,\,-x,\,y}.
\end{equation}
\end{pro}

\proof Let $n,k\in`N$. A composition of $n$ in $k$ parts is an ordered sequence $x_1,\dots,x_k$ of non negative integers, such that $x_1+\cdots+x_k=n$. The number of compositions of $n$ in $k$ parts (or of compositions of $n+k$ in $k$ positive integer parts) is
$\binom{n+k-1}{k-1}$.\par
We call \it run \rm of the chain $\S=(S(i))_{0\leq i\leq l}$, a maximal non-empty interval $I$ of $\cro{1,l}$ such that $(\Delta S_i)_{i\in I}$ is constant. The trajectories $\S$ of $T_{ab}^{j}(l,y)$ are composed by $j+\ind_{b=u}$ runs of $\u$ and $j+\ind_{a=d}$ runs of $\d$. The $\u$-runs form a compositions of $(l+y)/2$ (this is the number of steps $\u$) in positive integer parts, and the $\d$-runs form a composition of $(l-y)/2$ in positive integer parts. Hence, \[\#T_{ab}^{j}(l,y)=\binom{(l+y)/2-1}{j+\ind_{b=u}-1}\binom{(l-y)/2-1}{j+\ind_{a=d}-1},\] and Formula \eref{resum1} holds true.\par
The proofs of the other formulas are more tricky; the main reason for it is that the reflexion principle does not conserve the number of peaks. What still holds is, for any $x\geq 0$, $y\geq 0$, $j\geq 0$, $l\geq 0$
\[\#T_{ab}^{j,\geq }(l,x,y)=\#T_{ab}^{j}(l,x,y)-\#T_{ab}^{j,\vc}(l,x,y)\]
where $T_{ab}^{j,\vc}(l,x,y)$ is the set of trajectories belonging to $T_{ab}^{j}(l,x,y)$ that reach the level $-1$. Since $\#T_{ab}^{j}(l,x,y)=\#T_{ab}^{j}(l,0,y-x)$ is known, it remains to determine $T_{ab}^{j,\vc}(l,x,y)$. \\
We define two actions on the set of chains~:\\
\tb let $\S=(S(i))_{i\in\cro{0,l}}\in\W_l$. For any $t\in\cro{0,l}$
we denote by $\S'=\Ref(\S,t)$
the path $\S'=(S'_i)_{0\leq i\leq l}$ obtained from $\S$ by a reflexion from the abscissa $t$; formally~:
\[\left\{\begin{array}{ll}
S'(i)=S(i) &\textrm{ for any }0\leq i\leq t,\\
S'({i+t})=2S({i})-S({i+t})&\textrm{ for any }0\leq i\leq l-t\end{array}\right..\]
When $g$ is a function from $\W_l$ taking its values in $\cro{0,l}$, we write simply $\Ref(.,g)$ for the reflexion at abscissa $g(\S)$.\\
\tb let $\S=(S(i))_{i\in\cro{0,l}}\in\W_l$. For any $c$ and $d$ in $\mathbb{N}$, $0\leq c\leq d\leq l$,
we denote by $\S'=\Cont(\S,[c,d])$
the path $\S'=(S'(i))_{0\leq i\leq l-d}$ obtained from $\S$ by a contraction of the interval $[c,d]$~:
\[\left\{\begin{array}{ll}
S'(i)=S(i) &\textrm{ for any }0\leq i\leq c,\\
S'({c+i})=S(d+i)-S({d})+S(c)&\textrm{ for any }0\leq i\leq l-(d-c)\end{array}\right..\]
As before, we write $\Cont(.,[g_1,g_2])$ for the contraction of the interval $[g_1(\S),g_2(\S)]$.\\
\begin{figure}[htbp]
\centerline{\psfrag{c}{$c$}\psfrag{d}{$d$}\psfrag{t}{$t$}\psfrag{l}{$l$}\psfrag{l-(d-c)}{$l-(d-c)$}
\includegraphics[height=4.5cm]{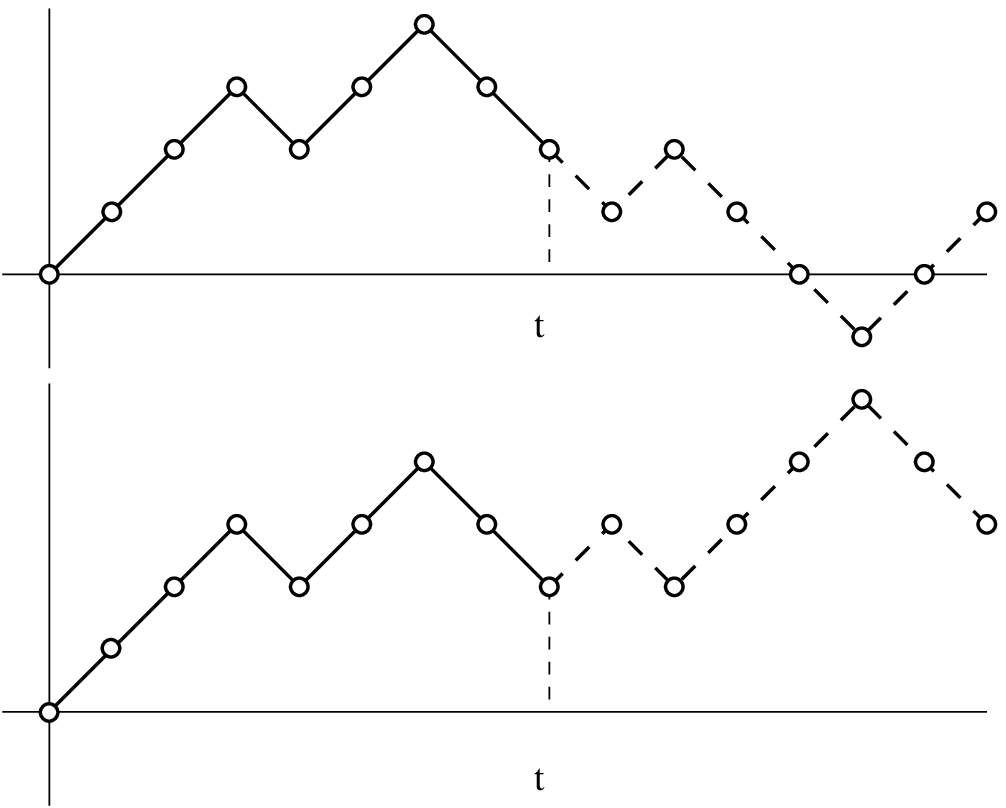}~~~~
\includegraphics[height=4.5cm]{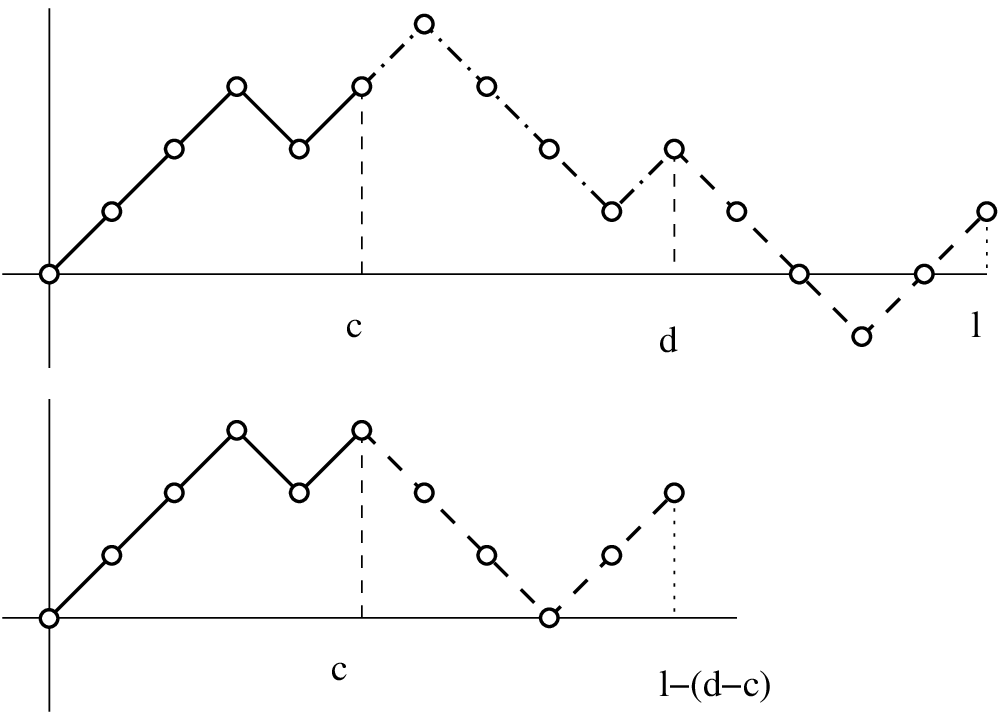} }
\captionn{On the first column $\S$ and $\Ref(\S,t)$, on the second column $\S$ and  $\Cont(\S,[c,d])$}
\end{figure}
We denote by ${\cal T}_{-1}(\S)=\inf\{j~:S(j)=-1\}$ the hitting time of $-1$ by $\S$. We proceed to a classification of the paths $\S$ from $\dis T_{ab}^{j,\vc}(l,x,y)$ according to the two first steps following ${\cal T}_{-1}(\S)$, that exist since $y$ is taken positive. We encode these two first steps following ${\cal T}_{-1}$ above the symbol $T$~: for any $\alpha,\beta\in\{\u,\d\}$, set
\[{\mathop{T}^{(\alpha\beta)}}_{ab}^{j,\vc}(l,x,y)=\{\S~: \S\in T_{ab}^{j,\vc}(l,x,y),~ \Delta S_{{\cal T}_{-1}(\S)+1}=\alpha,~\Delta S_{{\cal T}_{-1}(\S)+2}= \beta\}.\]
Hence, $T_{ab}^{j,\vc}(l,x,y)$ is the union of four elements of that type. Let us compute  $\dis\#{\mathop{T}^{(\alpha\beta)}}_{ab}$ for any $\alpha,\beta$. \par
For any $a\in\{\u,\d\}$, $x\geq 0$, $y\geq 1$, the following bijections (denoted by $\bij$) hold
\be
\tdu_{a\d}^{j,\vc}(l,x,y)\bij & \dis\tud_{a\u}^{j,\vc}(l,x,-2-y) &\bij T^{j-1}_{a\u}(l-2,x,-2-y)\\
\tdu_{a\u}^{j,\vc}(l,x,y)\bij & \dis\tud^{j+1,\vc}_{a\d}(l,x,-2-y)&\bij
T^{j}_{a\d}(l-2,x,-2-y).
\ee
On each line, the first bijection is $\Ref(.,{\cal T}_{-1})$, the second one is $\Cont(., [{\cal T}_{-1},{\cal T}_{-1}+2])$. Notice that this last operation does not create any peak because $\Delta S_{{\cal T}_{-1}(\S)}=\d$. The cardinalities of the sets in the right hand side are known, hence, $\dis\tdu_{\d\d},\tdu_{\u\d},\tdu_{\d\u},\tdu_{\u\u}$ are determined. Set 
$\dis{\mathop{T}^{(\alpha)}}_{ab}^{j,\vc}(l,x,y)=\dis{\mathop{T}^{(\alpha\u)}}_{ab}^{j,\vc}(l,x,y)\cup\dis{\mathop{T}^{(\alpha\d)}}_{ab}^{j,\vc}(l,x,y).$
We have
\be
\tu_{a\d}^{j,\vc}(l,x,y)\bij& \dis\td_{a\u}^{j-1,\vc}(l,x,-y-2)&\bij T_{a\u}^{j-1}(l-1,x,-y-1)\\
\tu_{a\u}^{j,\vc}(l,x,y)\bij& \dis\td_{a\d}^{j,\vc}(l,x,-y-2)&\bij T_{a\d}^{j}(l-1,x,-y-1)
\ee
The first bijection is $\Ref(.,{\cal T}_{-1})$, the second one $\Cont(., [{\cal T}_{-1},{\cal T}_{-1}+1])$.
Now,
\be
\#\tdd_{ab}^{j,\vc}&=& \#\tud_{ab}^{j+1,\vc}(l,x,y+2)\\
&=&\#\tu_{ab}^{j+1,\vc}(l,x,y+2)-\#\tuu_{ab}^{j+1,\vc}(l,x,y+2)\\
&=&\#\tu_{ab}^{j+1,\vc}(l,x,y+2)-\#\tdu_{ab}^{j+1,\vc}(l,x,y)
\ee
in the first line, we have replaced $\Delta S_{{\cal T}_{-1}}=\d$ by $\u$, in the third line we have replaced $\u\u$ by $\d\u$.~$\Box$

\subsection{Proof of Proposition \ref{card}}\label{ra} $(i)$ To build a path with $k$ peaks, dispose $k$ peaks, that is $k$ pairs $\u\d$.  Take a composition $x_1,\dots,x_{2(k+1)}$ of $n-2k$ in $2(k+1)$ non negative parts. Fill in now the  $k+1$ intervals between these peaks~: in the $l$th interval dispose $x_{2l-1}$ steps $\d$ and $x_{2l}$ steps $\u$.\\
$(ii)$ Assume $n=2N$ is even. To build a bridge with $k$ peaks, dispose $k$  pairs $\u\d$. Take two compositions $x_1,\dots,x_{(k+1)}$ and  $x_1',\dots,x_{(k+1)}'$  of $N-k$ in $k+1$ parts. Fill now the  $k+1$ intervals between these peaks~: in the $l$th interval dispose $x_l$ steps $\d$ and  $x_l'$ steps $\u$.
\[(iii)~~~~~~~~~\#{\cal E}_n^{(k)}=\#T_{\u\d}^{k,\,\geq}(n-1,\,0,\,1)+\#T_{\u\u}^{k-1,\,\geq}(n-1,\,0,\,1)~~~~~~~~~~~~~~~ \]
For $(iv)$, one may use the bijections described in Section \ref {bml}, or proceed to a direct computation as follows; first,
\ben
\label{susu} \#{\cal M}_n^{(k)}&=&\#{\cal E}_n^{(k)}+\sum_{y\geq 1}T_{\u\u}^{k,\geq}(n,0,y)+T_{\u\d}^{k,\geq}(n,0,y).
\een
Denote by $W(n,k)$ the sum in \eref{susu}. The integer $W(n,k)$ is the number of meanders with length $n$, ending in a positive position. Using that $\binom nk+\binom n{k-1}=\binom{n+1}k$, we have 
\be
W(n,k)&=&\sum_{y\geq 1}  \g{n,k,k-1,0,y}-\g{n,k-1,k,0,y}+\g{n,k-1,k-1,0,y}-\g{n,k-2,k,0,y}.\\
&=& \sum_{y\geq 1}\binom{\dis\frac{n+y}{2}}{k}\binom{\dis\frac{n-y}{2}-1}{k-1}-\binom{\dis\frac{n+y}{2}}{k-1}\binom{\dis\frac{n-y}{2}-1}{k}.
\ee
Let  $a,b,c,k$ be positive integers, the following formula holds
\begin{equation}\label{su2}
\sum_{y\geq c} \binom{a+y}{k}\binom{b-y}{k-1}-\sum_{y\geq c}\binom{a+y}{k-1}\binom{b-y}{k}=\binom{a+c}{k}\binom{b-c+1}{k}.
\end{equation}
Indeed: a term in the first sum counts the number of ways to choose $2k$ items among $a+b+1$, choosing, the $a+y+1$th, and $k$ items among the $a+y$ first ones, when,  a term in the second sum counts the number of ways to choose $2k$ items among $a+b+1$, choosing, the $a+y+1$th, and $k-1$ items among the $a+y$ first ones. 
The choices counted by the first sum but not by the second one, are those where exactly $k$ items are chosen among the $a+c$ first ones.\par
We need to consider the two cases $n$ even and $n$ odd~:\\
$\bullet$ If $n$ is even, using that $S(n)=n \mod 2=y\mod 2$, set $N=n/2,$ $Y=y/2$ in the sum, 
\[W(2N,k)=\sum_{Y\geq 1} \binom{N+Y}{k}\binom{N-Y-1}{k-1}-\binom{N+Y}{k-1}\binom{N-Y-1}{k}=\binom{N+1}{k}\binom{N-1}{k}.\]
If $n=2N+1$, take $y=2Y+1$ in the sum 
\[W(2N+1,k)=\sum_{Y\geq 0} \binom{N+Y+1}{k}\binom{N-Y-1}{k-1}-\binom{N+Y+1}{k-1}\binom{N-Y-1}{k}=\binom{N+1}k\binom Nk.~~\Box\]

\section{Asymptotic considerations and proofs}
\label{ACP}
We first recall a classical result of probability theory, simple consequence of Billingsley \cite[Theorem 7.8]{BIL},
that allows to prove the weak convergence in $`R^d$ using a local limit theorem.
\begin{pro} \label{bill}
Let $k$ be a positive integer and for each $i\in\cro{1,k}$,   $(\alpha_n^{(i)})$ a sequence of real numbers such that $\alpha_n^{(i)}\sous{\tend}{n} +\infty$.
For any $n$, let $X_n=(X_n^{(1)},\dots, X_n^{(k)})$ be a $\mathbb{Z}^k$-valued random variable. If $\forall \lambda=(\lambda_1,\dots,\lambda_k)\in`R^k$,
\[\alpha_n^{(1)}\dots\alpha_n^{(k)}`P\big((X_n^{(1)},\dots,X_n^{(k)})=(\floor{\lambda_1\alpha_n^{(1)}},\dots,\floor{\lambda_k\alpha_n^{(k)}})\big)\tend \phi(\lambda_1,\dots,\lambda_k)\]
where $\phi$ is a density probability on $`R^k$, then $(X_n^{(1)}/\alpha_n^{(1)},\dots,X_n^{(k)}/\alpha_n^{(k)})\dd X$ where  $X$ is a random variable with density $\phi$.
\end{pro}
Proposition \ref{CVD} is a simple consequence of this proposition, since by \eref{card2} and Proposition \ref{card},  the application of the Stirling formula simply yields
\[\sqrt{n}\,`P_n^{\,\x}(\floor{\Lambda_n- n/4}=\floor{t\sqrt{n}})\xrightarrow[~n~]{} \frac{1}{\sqrt{2\pi/16}}\exp(-8t^2),\]  for any $\x\in\{\w,\b,\e,\m\}$ and any $t\in`R$. Note that under $`P_n^{\w}$, one may also compute the limiting distribution using that $\Lambda_n(\S)=\sum_{i=1}^{n-1} \ind_{i\in \S_{\wedge}}$, which is a sum of Bernoulli random variables with an easy to handle dependence. 

\subsection{Finite dimensional distribution of the Brownian processes}
\label{fdd}
\bf Notation \rm For any sequence $(o_i)_i$ indexed by integers, the sequence $(\Delta o_i)$ is defined by $\Delta o_i=o_i-o_{i-1}$ and  $(\bar\Delta o_i)$ by $\bar \Delta o_i=o_i+o_{i-1}$. \par
For any $t>0$ and $x,y\in`R$, set
\[\p_t(x,y)=\frac{1}{\sqrt{2\pi t}}\exp\l(-\frac{(y-x)^2}{2t}\r).\]
Let $\ell\geq 1$ and let $(t_1,\dots, t_\ell) \in[0,1]^\ell$ satisfying $0< t_1<\dots <t_{\ell-1}<t_\ell:=1$. The distributions of the a.s. continuous processes $\w,\b,\e,\m$ are characterized by their finite dimensional distributions. Let $f^{\x}_{t_1,\dots,t_k}$ be the density of $(\x_{t_1},\dots,\x_{t_k})$ with respect to the Lebesgue measure on $`R^k$.
We have
\be
f^{\w}_{t_1,\dots, t_\ell}(x_1,\dots, x_\ell)&=& \prod_{i=1}^\ell \p_{\Delta t_i}(x_{i-1},x_i), ~~~\textrm{ with }x_0=0 \textrm{ by convention}, \\
f^{\b}_{t_1,\dots, t_{\ell-1}}(x_1,\dots, x_{\ell-1})&=& \sqrt{2\pi}\,f^{\w}_{t_1,\dots, t_{\ell-1}}(x_1,\dots, x_{\ell-1})\p_{\Delta t_\ell}(x_{\ell-1},0),\\
f^{\m}_{t_1,\dots, t_\ell}(x_1,\dots, x_\ell)&=& \sqrt{2\pi}\frac{x_1}{t_1}\p_{t_1}(0,x_1) \l(\prod_{i=2}^{\ell} \p_{\Delta t_i}(x_{i-1},x_i)-\p_{\Delta t_i}(x_{i-1},-x_i)\r)\ind_{x_1,\dots,x_\ell\geq 0},\\
f^{\e}_{t_1,\dots, t_{\ell-1}}(x_1,\dots, x_{\ell-1})&=&f^{\m}_{t_1,\dots, t_{\ell-1}}(x_1,\dots, x_{\ell-1})\frac{x_{\ell-1}}{\Delta t_{\ell}}\p_{\Delta t_{\ell}}(x_{\ell-1},0).
\ee
We end this section with two classical facts: first $\p_t(x,y)=\p_t(0,y-x)$ and, for any $\alpha>0$, $f^{\alpha\x}_{t_1,\dots,t_k}=\alpha^{-k} f^{\x}(\alpha^{-1} t_1,\dots,\alpha^{-1} t_k)$.

\subsection{Finite dimensional convergence}
We will show that for any $\x\in\{\w,\b,\e,\m\}$, under $`P_n^{\x}$, for any $\ell \in\{1,2,3,\dots\}$ and $0<t_1<\dots<t_{\ell-1}<t_\ell:=1$
\[\l(s_n(t_1),\dots,s_n(t_\ell),\lambda_n(t_1),\dots,\lambda_n(t_\ell)\r)\dd \l(\x(t_1),\dots,\x(t_\ell),\widehat{\b}(t_1),\dots,\widehat{\b}(t_\ell)\r),\] 
where $t_\ell:=1$ has been chosen for computation convenience, and $(\x,\widehat{\b})$ has the prescribed distribution (as in Theorem \ref{pri}). This implies the convergence of the finite dimensional distribution  in Theorem \ref{pri}. \\
In order to handle easily the binomial coefficients appearing in $\g{l,j_1,j_2,x,y}$ that involve half integers, we proceed as follows. Let  $N=\floor{n/2}$ and let ${\cal E}_n=n-2N=n \mod 2$. 
For $i\in\cro{1,\ell-1}$, let $t_i^{(n)}$ be defined by
\[2Nt_i^{(n)}:= 2\floor{nt_i/2},\]
 $t_{0}^{(n)}=0$, and $t_\ell^{(n)}$ by $2Nt_\ell^{(n)}=2N+{\cal E}_n=n$ (notice that $t_\ell^{(n)}$ is in $\{1,1+1/n\}$).\par
 Using that for any $i$, $|2Nt_i^{(n)}-nt_i|\leq 2$, 
we have clearly under $`P_n^{\x,(\k{n})}$,
\begin{equation}
\frac{S(2Nt_i^{(n)})-S(nt_i)}{\beta_n}\sur{\tend}{proba.}0 \textrm{ and }\frac{\Lambda_{2Nt_i^{(n)}}-\Lambda_{nt_i}}{\gamma_n}\sur{\tend}{proba.}0\end{equation}
since $\gamma_n$ and $\beta_n$ goes to $+\infty$.
From now on, we focus on the values of the processes on the discretization  points $2Nt_i^{(n)}$.
For any $i\in\cro{1,\ell-1}$, set
\[\widetilde{\Lambda}_{i}:=\# \S_{\wedge}\cap\cro{2Nt_{i-1}^{(n)}+1,2Nt_{i}^{(n)}-1}\]
the number of peaks lying strictly between $2Nt_{i-1}^{(n)}$ and $2Nt_{i}^{(n)}$.\par 
In order to obtain a local limit theorem, we are interested by the number of trajectories passing via some prescribed positions. 

\subsubsection{Case  $\x=\w$} \rm 
Let $0=u_0, u_1,\dots,u_\ell,v_1,\dots,v_{\ell-1}$ be fixed real numbers. 
Set \[\Theta:=(t_1,\dots,t_{\ell-1},u_1,\dots,u_\ell,v_1,\dots,v_{\ell-1})\]
and for any $i\in\cro{1,l-1}$, set 
\[\k{n}^{(i)}=\floor{\k{n}\Delta t_{i}}+\floor{v_{i}\gamma_n},\]
and
\be
A_n^\w(\Theta)&:=& \{\S\in \W_n^{(\k{n})}, S(2Nt_i^{n})=2\floor{u_i\beta_n},\textrm{ for }i\in\cro{1,\ell-1}, 
S(n)=2[u_\ell\beta_n]+{\cal E}_n,\\
&&\widetilde{\Lambda}_{i}= \k{n}^{(i)}\textrm{ for }i\in\cro{1,\ell-1}\}.
\ee
For any $i\in\cro{1,l-1}$,  denote by $a_i(\S)=\Delta S_{2Nt_{i-1}^{n}}+1$ and $b_i(\S)=\Delta S_{2Nt_i^{n}}$ the first and last increments of the $i$th part of $\S$ between the discretization points.
Some peaks may appear in the positions $2Nt_i^{n}$, and then, we must take into account the pairs $(b_i,a_{i+1})$ to compute the cardinality of $A_n^\w(\Theta)$. For any $\S\in A_n^\w(\Theta)$, the number of peaks in $\cro{2Nt_{\ell-1}+1,n-1}$ is 
\[\k{n}^{(\ell)}=\k{n}-\sum_{i=1}^{\ell-1} \k{n}^{(i)}-\sum_{i=1}^{\ell-1}\indic_{(b_{i},a_{i+1})=(\u,\d)}.\]
We have $\#A_n^\w(\Theta)=$
\begin{equation} \label{grosse}
\sum_{c} \l(\prod_{i=1}^{\ell-1}\#T_{a_i,b_i}^{\k{n}^{(i)}}(2N \Delta t_{i}^{n},2\floor{u_{i-1}\beta_n},2\floor{u_i\beta_n})\r)\#T_{a_\ell,b_\ell}^{\k{n}^{(l)}}(2N\Delta t_{\ell}^{n},2[u_{\ell-1}\beta_n],2[u_\ell\beta_n]+{\cal E}_n)
\end{equation}
where the sum is taken over every $c:=((a_1,b_1),\dots,(a_\ell,b_\ell))\in\{u,d\}^\ell$.

In order to evaluate the sum \eref{grosse}, we introduce some binomial random variables $B(l,p_n)$ with parameters $l$ and $p_n$
and we use the following version of the local limit theorem
\begin{lem}\label{lll}
 Let $(l(n))$ be a sequence of integers going to $+\infty$ and $\sigma^2_n=l(n)p_n(1-p_n)$. We have
\begin{equation}
`P(B(l(n),p_n)=m)= \binom{l(n)}{m}p_n^{m}(1-p_n)^{l(n)-m}=\frac{1}{\sigma_n\sqrt{2\pi}}\exp\l(-\frac{(m- l(n) p_n)^2}{2 \sigma_n^2}\r)+o(1/\sigma_n).\end{equation}\end{lem}
This may be proved using Stirling formula.
As a consequence, if $(a_n)$, $(a_n')$, $(a_n'')$ are sequences of integers such that
\[a_n-nt/2=O(1) \textrm{ for }t\in(0,1),~~~ \frac{a_n'}{\beta_n}\tend u, ~~~  \frac{a_n''-t \k{n}}{\gamma_n}\tend v,\]then
\[`P(B(a_n+a_n',p_n)=a_n'')
\sim\frac{1}{\sqrt{\pi t}\gamma_n}\exp\l(-\frac{(v-u)^2}{t}\r).\]
We then get easily that
\[\#T_{a_i,b_i}^{\k{n}^{(i)}}(2N \Delta t_{i}^{n},2\floor{u_{i-1}\beta_n},2\floor{u_i\beta_n})\,(1-p_n)^{2N\Delta t_i^n-2}\l(\frac{p_n}{1-p_n}\r)^{2\k{n}^{(i)}-\ind_{b_i=\d}-\ind_{a_i=\d}}\]
is equivalent to $\frac{1}{\pi\Delta t_i \gamma_n ^2}\exp\l(-2\frac{v_i^2+\Delta u_i^2}{\Delta t_i}\right)=2\gamma_n^{-2}\p_{\Delta t_i}(0,2v_i)\p_{\Delta t_i}(0,2\Delta u_i)$ and
\[\#T_{a_\ell,b_\ell}^{\k{n}^{(\ell)}}(2N\Delta t_{\ell}^{n},2[u_{\ell-1}\beta_n],2[u_\ell\beta_n]+{\cal E}_n)\,(1-p_n)^{n-2Nt_{\ell-1}^{(n)}-2}\l(\frac{p_n}{1-p_n}\r)^{2\k{n}^{(\ell)}-\ind_{b_\ell=\d}-\ind_{a_\ell=\d}}\]
is equivalent to  $\frac{1}{\pi \Delta t_\ell \gamma_n ^2}\exp\l(-2\frac{(\sum_{i=1}^{\ell-1} v_i)^2+\Delta u_\ell^2}{\Delta t_\ell}\right)=2\gamma_n^{-2}\p_{\Delta t_i}(2\sum_{i=1}^{\ell-1} v_i,0)\p_{\Delta t_i}(0,2\Delta u_\ell)$.

Since $p_n^{2\k{n}+1}(1-p_n)^{n-2\k{n}}\,\,\# \W_n^{(\k{n})}
\sim 1/(\gamma_n\sqrt{2\pi})$
we obtain
\[`P^{\w,(\k{n})}_n(A_n^\w(\Theta))=\frac{\#A_n^\w(\Theta)}{\# \W_n^{(\k{n})}}\sim 
c_n^\w f^{\w/2}_{t_1,\dots,t_\ell}(u_1,\dots,u_l)f^{\b/2}_{t_1,\dots,t_{\ell-1}}(v_1,v_1+v_2,\dots,v_1+\dots+v_{\ell-1})\]
where
\be
c_n^\w&:=& 2^{1-\ell}\gamma_n^{1-2\ell}p_n(1-p_n)^{2\ell}\sum_{c}  \l(\frac{p_n}{1-p_n}\r)^{\ind_{a_1=u}+\ind_{b_\ell=\d}+\sum_{i=1}^{\ell-1} 2\ind_{(b_i,a_{i+1})=(u,d)}+\ind_{b_i=d}+\ind_{a_{i+1}=u}}\\
&=& \gamma_n^{1-2\ell}p_n^\ell=\gamma_n^{1-\ell}\beta_n^{-\ell}
\ee
The contribution of the sum over $c$ has been computed as follows:
\be
\sum_c x^{f(c)}&=&\l(\sum_{(a_1,b_\ell)\in \{+1,-1\}^2}x^{\ind_{a_1=u}+\ind_{b_\ell=\d}}\r)\prod_{i=1}^{\ell-1}\sum_{(b_i,a_{i+1})\in \{+1,-1\}^2} x^{\ind_{(b_i,a_{i+1})=(u,d)}+\ind_{b_i=d}+\ind_{a_{i+1}=u}}\\
&=& (1+2x+x^2)(2x+2x^2)^{\ell-1} =2^{\ell-1}x^{\ell-1}(1+x)^{\ell+1}.
\ee
Finally, this says that 
\[\beta_n^\ell \gamma_n^{\ell-1}`P^{\w,(\k{n})}_n\l(\frac{S(2Nt_i)}{2}=\floor{u_i\beta_n}, i\in\cro{1,\ell},
\widetilde{\Lambda}_i-\floor{\k{n}\Delta t_i}=\floor{v_i\gamma_n}, i\in\cro{1,\ell-1}\r)\to\]\[ f^{\w/2}_{t_1,\dots,t_\ell}(u_1,\dots,u_l)f^{\b/2}_{t_1,\dots,t_{\ell-1}}(v_1,v_1+v_2,\dots,v_1+\dots+v_{\ell-1}).\]
Hence by Proposition \ref{bill} and \eref{sl}, and taking into account that for any $i$, $\lambda_n(t_i)-2\sum_{j=1}^i\frac{\widetilde{\Lambda}_i-\floor{\k{n}\Delta_{t_i}}}{\gamma_n}\sur{\tend}{proba}0$, this allows to conclude to the finite dimensional convergence in Theorem \ref{pri} in the case $\x=\w$.
\begin{comm}To compute a local limit theorem under the other distributions the numbers $u_1,\dots,u_\ell,v_1,\dots,v_\ell$ and the set $A_n^{\w}$ have to be suitably change. First, in each case, the set ${\W}_n^{(\k{n})}$ has to be replaced by the right set. \\
$\bullet$ In the case of excursions and bridges, $n$ is an even number and $u_\ell$ is taken equal to 0. \\
$\bullet$ In the case of excursions  $a_1=\u$, $b_\ell=\d$\\
$\bullet$ In the case of excursions and meanders all the reals $u_i$ are chosen  positive. Moreover, $T^{\geq}$ must replace $T$ in the summation \eref{grosse}.\par
Up to these changes, the computation are very similar to the case of Bernoulli chains.
\end{comm}

\subsubsection{Case $\x=\b$} The computation is very similar to the previous case; the only differences are~: here $n=2N$ is even, $\# \B_n^{(\k{n})}=\binom{N}{\k{n}}^2$, we set  $u_\ell$ to 0 and we take $\Theta':=(t_1,\dots,t_{\ell-1},u_1,\dots,u_{\ell-1},1,v_1,\dots,v_{\ell-1}).$ We get
\be
`P^{\b,(\k{n})}_n(A_n^\b(\Theta'))\sim c_n^\b  f^{\b/2}_{t_1,\dots,t_{\ell-1}}(u_1,\dots,u_{\ell-1})f^{\b/2}_{t_1,\dots,t_{\ell-1}}(v_1,v_1+v_2,\dots,v_1+\dots+v_{\ell-1}),
\ee
where $c_{n}^\b:= c_{n}^\w \gamma_n p_n^{-1}=\gamma_n^{1-\ell} \beta_{n}^{1-\ell}$.
\subsubsection{Case $\x=\e$}  In this case $n=2N$ is even, $u_\ell=0$, $a_1=\u,b_\ell=\d$. In order to avoid some problems with the formulas provided in Proposition \ref{denu} that have to be handled with precautions when $x$ or $y$ are 0, we will compute the local limit theorem ``far'' from 0. This will however suffice to conclude. For $i\in\cro{1,\ell-1}$ we take $u_i>0$, and $\beta_n$ is assumed large enough so that $\floor{u_i \beta_n}>0$. For the calculus in the case of $\x=\e$, in formula \eref{grosse}, we replace $T$ by $T^{\geq}$. Finally, $\#{\cal E}_n^{(\k{n})}=\frac1N\binom N{\k{n}}\binom N{\k{n}-1}$. \par

We first treat the contribution of the non extreme parts of the trajectories, namely, $i\in\cro{2,l-1}$, 
\[\#T_{a_i,b_i}^{\k{n}^{(i)},\geq }(2N \Delta t_{i}^{n},2\floor{u_{i-1}\beta_n},2\floor{u_i\beta_n}))(1-p_n)^{2N\Delta t_i^n-2}\l(\frac{p_n}{1-p_n}\r)^{2\k{n}^{(i)}-\ind_{b_i=\d}-\ind_{a_i=\u}}\]
is equivalent to $2\gamma_n^{-2}\p_{\Delta t_i}(0,2v_i)\l(\p_{\Delta t_i}(0,2\Delta u_i)-\p_{\Delta t_i}(0,2\bar{\Delta}u_i)\r)$.

Let us consider $i=1$. Notice that $\#T_{\u\u}^{j,\geq}(l,0,y)$ and $\#T_{\u\d}^{j,\geq}(l,0,y)$ may be very different~:
\[T_{\u\u}^{j,\geq}(l,0,y)=\binom{\frac{l+y}2-1}{j-1}\binom{\frac{l-y}2-1}{j-1}\frac{y}{j}~\textrm{ and }~T_{\u\d}^{j,\geq}(l,0,y)=\binom{\frac{l+y}2-1}{j-1}\binom{\frac{l-y}2-1}{j-1}\frac{2jy+l-y}{j(l+y-2j+2)}\]
\[A_{b_1}(n):=\#T_{\u,b_1}^{\k{n}^{(1)},\geq }(2N \Delta t_{1}^{n},0,2\floor{u_1\beta_n}))(1-p_n)^{2N\Delta t_1^n-2}\l(\frac{p_n}{1-p_n}\r)^{2\k{n}^{(1)}-\ind_{b_1=\d}}\]
We notice that under $(H)$ 
\[n=o(\k{n} \beta_n),~~ \gamma_n=o(n-2\k{n}),~~ \gamma_n\leq \beta_n,~~  \gamma_n=o(\k{n}).\]
We then get,
\[A_{b_1}(n)\sim \frac{1}{n\gamma_n(1-p_n)^2}\frac{8u_1}{t_1}\p_{\Delta t_1}(0,2u_1)\p_{\Delta t_1}(0,2v_1).\]
The case $i=\ell$ is treated with the same method. 
We obtain \[`P^{\e,(\k{n})}_n(A_n^\e(\Theta'))=\frac{\#A_n^\e(\Theta')}{\# \E_n^{(\k{n})}}
\sim c_n^\e f^{\e/2}_{t_1,\dots,t_{\ell-1}}(u_1,\dots,u_{\ell-1})f^{\b/2}_{t_1,\dots,t_{\ell-1}}(v_1,v_1+v_2,\dots,v_1+\dots+v_{\ell-1})\]
where $c_n^\e=\gamma_n^{-\ell+1}\beta_n^{-\ell+1}$.
\subsubsection{Case $\x=\m$} The computation is the same as in the case $\x=\e$, except that $u_\ell$ is taken $>0$ (and $n$ large enough such that $\floor{u_\ell\beta_n}\geq 1$)· The last piece in the decomposition of meanders is of the same type as a standard excursion piece.
We obtain \[`P^{\m,(\k{n})}_n(A_n^\m(\Theta))=\frac{\#A_n^\m(\Theta')}{\# \M_n^{(\k{n})}}
\sim c_n^\m f^{\m/2}_{t_1,\dots,t_{\ell-1}}(u_1,\dots,u_{\ell-1})f^{\b/2}_{t_1,\dots,t_{\ell-1}}(v_1,v_1+v_2,\dots,v_1+\dots+v_{\ell-1})\]
where $c_n^\m=\gamma_n^{-\ell+1}\beta_n^{-\ell}$. 

\section{Tightness}

We begin with some recalls of some classical facts regarding tightness in $C[0,1]$. First, tightness and relative compactness are equivalent in $C[0,1]$ (and in any Polish space, by Prohorov). 
Consider the function modulus of continuity, 
\[\app{\omega}{[0,1]\times C[0,1]}{`R^+}{(\delta,f)}{\omega_{\delta}(f)}\]
defined by
\[\omega_{\delta}(f):=\sup_{s,t\in[0,1], |s-t|\leq \delta} |f(t)-f(s)|.\]
A sequence of processes $(x_n)$ such that $x_n(0)=0$ is tight in $C[0,1]$,  if for any $`e>0, \eta>0$, there exists $\delta>0$, such that for $n$ large enough
\[`P( \omega_{\delta}(x_n)\geq `e)\leq \eta.\]
If the sequences $(x_n)$ and $(y_n)$ are tight in $C[0,1]$ (and if for each $n$, $x_n$ and $y_n$ are defined on the same probability space $\Omega_n$), then the sequence $(x_n,y_n)$ is tight in $C([0,1])^2$.
We will use this result here, and prove the tightness separately for $(s_n)$ and $(\lambda_n)$ for every model $`P^\x_n$.\par
We say that a sequence $(x_n)$ in $C[0,1]$ is tight on $[0,1/2]$ if the sequence of restrictions $(x_n|_{[0,1/2]})$ is tight in $C[0,1/2]$. We would like to stress on the fact that we deal only with processes piecewise interpolated (on intervals $[k/n, (k+1)/n]$); for these processes,  for $n$ large enough such that $1/n<\delta$, 
\[\Big\{\sup\{ |x_n(t)-x_n(s)|,s,t\in\{k/n, k\in\cro{0,n}\}, |s-t|\leq \delta\}\leq `e/3\Big\}\imp \{\omega_{\delta}(x_n)\leq `e\}. \]
In other words, one may assume that $s$ and $t$ are discretization points, in our proofs of tightness.

We recall a result by Petrov \cite[Exercise 2.6.11]{PET}~:
\begin{lem}\label{berbound}
 Let $(X_i)_i$ be i.i.d. centered random variables, such that $`E(e^{tX_1})\leq e^{g t^2/2}$ for $|t|\in[0,T]$ and  $g>0$. Let $Z_k=X_1+\dots+X_k$. Then 
\[`P(\max_{1\leq k\leq N} |Z_k|\geq x)\leq 2\left\{
\begin{array} {ll}
\exp(-x^2/2Ng)&\textrm{ for any } x\in[0,NgT]\\
\exp(-Tx)&\textrm{ for any } x\geq NgT.
\end{array}\right.
\]
\end{lem}

The tightness of $(s_n,\lambda_n)$ is proved as follows: first, under $`P_n^{\w,(\k{n})}$, the passage via an alternative model of ``simple random walk'' allows to remove the conditioning by $\Lambda_n=\k{n}$. Then, the tightness under $`P_n^{\b,(\k{n})}$ is deduced from that under $`P_n^{\w,(\k{n})}$, thanks to the fact, that the conditioning by $S(n)=0$ does not really change the distribution the first half of the trajectories. The tightness under $`P_n^{\e,(\k{n})}$ and $`P_n^{\m,(\k{n})}$ are then obtained from that under $`P_n^{\b,(\k{n})}$, by some usual trajectory transformations that preserve the main properties of the variations and  peak distributions of the trajectories.

\subsection{A correspondence between simple chains and Bernoulli chains}

We denote by $\HH_{n}$ the set of ``simple chains'', starting from 0 and having $n+1$ steps~: 
\[\HH_{n}=\{\H=(H_i)_{0\le i \le n+1}~: H_0=0, H_{i+1}=H_i \textrm{ or } H_{i+1}=H_i+1 \textrm{ for any }i\in\cro{0,n}\}.\]

\begin{figure}[htbp]
\psfrag{0}{0}\psfrag{S}{$\S$}\psfrag{H}{$\H$}\psfrag{2}{2}\psfrag{4}{4}\psfrag{6}{6}\psfrag{7}{7}
\psfrag{n+1}{$n+1$}
\centerline{\includegraphics[height=5.6cm]{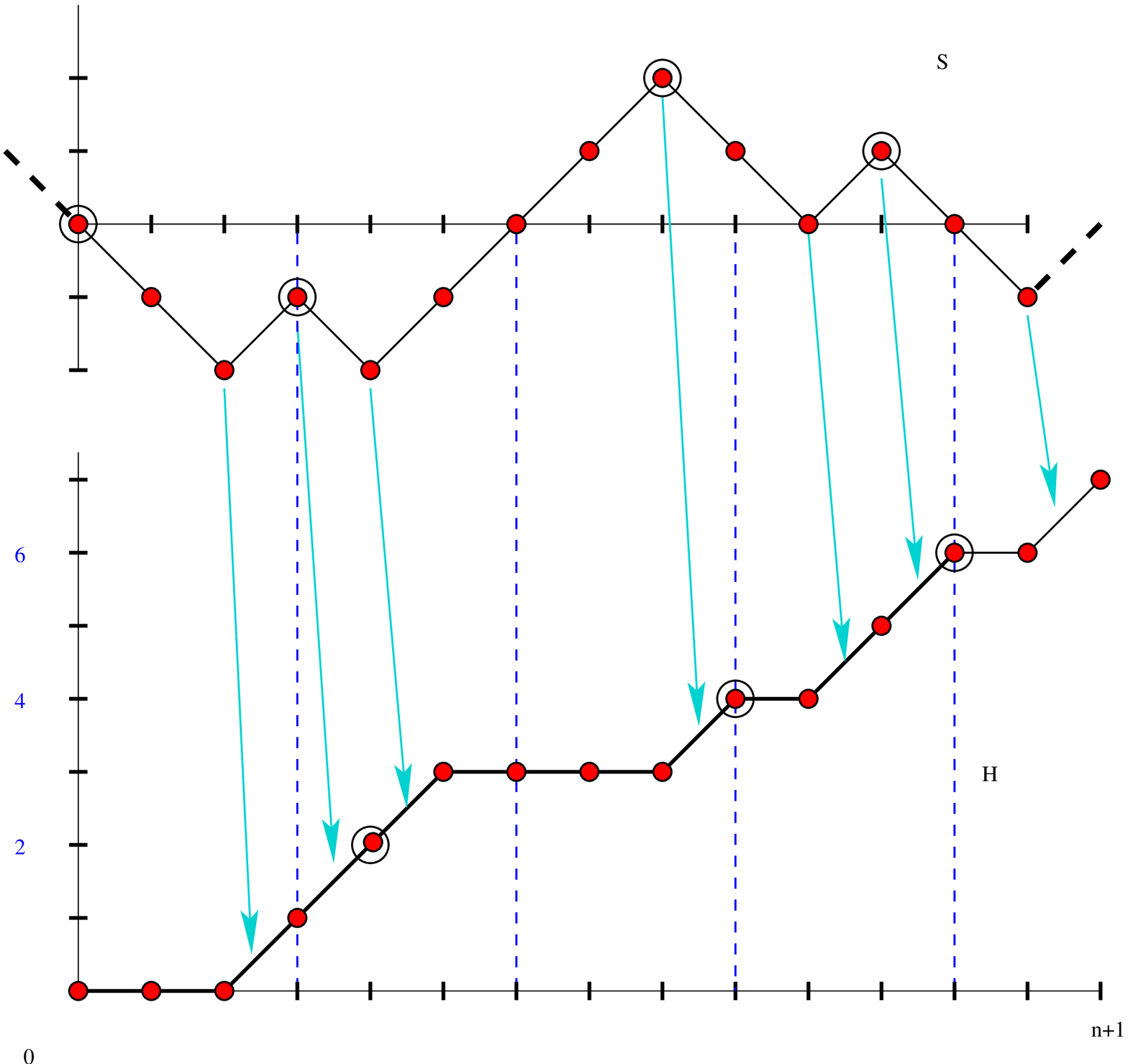}}
\caption{\label{cor}Correspondence between simple chains and Bernoulli chains}
\end{figure}

We consider the application 
\[\app{\Phi_n}{\W_n}{\HH_n}{\S}{\H=\Phi_n(\S)}\]
where $\H$ is the simple chain with increments: for any $i\in\cro{1,n+1}$,
\[\left\{
\begin{array}{ll}
\textrm{if }\Delta S_i\neq \Delta S_{i-1} &\textrm{ then } \Delta H_i=1\\
\textrm{if }\Delta S_i= \Delta S_{i-1} &\textrm{ then } \Delta H_i=0
\end{array}
\right.
\]
where by convention $\Delta S_{0}=-1$ and $\Delta S_{n+1}=1$ (see illustration on Figure \ref{cor}).

The mapping $\Phi_n$ is a combinatorial trick. Obviously, the application $\S\mapsto \H$ where $\H$ is defined by $\Delta H_i=(\Delta S_i+1)/2$ is a bijection from ${{\cal W}_n}$ onto ${{\cal H}_{n-1}}$. The application $\Phi_n$ is then certainly not a bijection (it is an injection). But, $\Phi_n$ owns some interesting properties that will really simplify our task.\par
Each increasing step in $\H$ corresponds to a changing of direction in $\S$. Since  $\Delta S_{0}=-1$, the first one corresponds then to a valley, and the last one to a peak (which can not be in position $n$, since $\Delta S_{n+1}=1$). 
Hence, for any $j\in\cro{0,n}$,
\be\label{are}
\Lambda_j(\S)=\#\S_{\wedge}\cap\cro{0,j}=\l\lfloor\frac{H_{j+1}}2\r\rfloor. 
\ee
Hence, 
\[\S_{\wedge}=\{{\cal T}_{2l}(\H)-1, l\in\mathbb{N}\} ~~\textrm{ and }~~(-\S)_{\wedge}=\{{\cal T}_{2l+1}(\H)-1, l\in\mathbb{N}\},\]
where  $(-\S)_{\wedge}$ is the set of valleys of $\S$ and where ${\cal T}_l(\H)=\inf\{j, H_j=l\}$ is the hitting time by $\H$ of the level $l$.
The process $\S$ may then be described with $\H$: 
\begin{equation}\label{gege}
S(k)=\sum_{i=0}^{H_k-1} (-1)^{i+1} \l({\cal T}_{i+1}(\H)- {\cal T}_{i}(\H)\r)+(-1)^{H_k+1}(k-{\cal T}_{H_k}(\H)).
\end{equation}
To end these considerations, consider now the subset of simple chains with $k$ increasing steps, 
\[\HH_n^{k}=\{\H\in \HH_n, H_{n+1}=k\},\]
and focus on $\HH_n^{2k+1}$. Each element $\H\in \HH_n^{2k+1}$ is image by  $\Phi_n$ of a unique trajectory $\S$ that has $k$ peaks  in $\cro{1,n-1}$ and $k+1$ valleys in $\cro{0,n}$ (that may be in position 0 and $n$ by construction), in other words, to a trajectory of $\W_n^{(k)}$.\par
This may alternatively be viewed as follows: to build a trajectory of $\W_n^{(k)}$ choose $2k+1$ integers $i_1<i_2<\dots<i_{2k+1}$ in the set $\cro{0,n}$. Then construct a trajectory from $\W_n$ in placing a valley in $i_{1},i_3,\dots,i_{2k+1}$, a peak in $i_2,i_4,\dots,i_{2k}$ and fill in the gaps between these points by straight lines.  Hence
\begin{lem}\label{restr} For any $k\in\cro{0,\floor{n/2}}$, the restriction of $\Phi_n$ on $\W_n^{(k)}$ is a bijection onto ${\cal H}_n^{2k+1}$.
\end{lem}
For any $p\in [0,1]$, let $`Q^{n}_p$ be the distribution on ${\cal H}_{n}$ of the Bernoulli random walks with $n+1$ i.i.d. increments, Bernoulli $B(p)$ distributed (that is $`Q^{n}_p(\Delta H_i=1)=1-`Q^{n}_p(\Delta H_i=0)=p$). For any $\H$ in ${\cal H}_n$, 
\[\mathbb{Q}^{n}_{p}(\{\H\})=p^{H_{n+1}}(1-p)^{n+1-H_{n+1}},\]
and then, $\mathbb{Q}^{n}_{p}$ gives the same weight to the trajectories ending at the same level.  Hence the conditional distribution $`Q_p^{n}(~.~| {\cal H}_n^{2k+1})$ is the uniform law on ${\cal H}_n^{2k+1}$. On the other hand, since $\mathbb{`P}_n^{\w,(k)}$ is the uniform distribution on $\W^{(k)}_n$, by Lemma \ref{restr}, $\mathbb{`P}_n^{\w,(k)}\circ \Phi_n^{-1}$ is also the uniform law on ${\cal H}_n^{2k+1}$. Hence
\begin{lem}   For any $p\in(0,1)$,  $n\in\mathbb{N}$, $k\in\cro{0,\floor{n/2}}$,  $`Q_p^{n}(~.~| {\cal H}_n^{2k+1})=\mathbb{`P}_n^{\w,(k)}\circ \Phi_n^{-1}$.
\end{lem}
Using simple properties of binomial distribution, the value of $p$ that maximizes $`Q_p^{n}({\cal H}_n^{2\k{n}+1})$ is $\tilde{p}_n=(2\k{n}+1)/(n+1)$.  This morally explains why in Section \ref{ACP}, $p_n$ appears as a suitable parameter. For sake of simplicity, we will work again with $p_n=2\k{n}/n$ instead of $\tilde{p}_n$. We will see that under $`Q_{p_n}^{n}$, 
 the conditioning by ${\cal H}_n^{2k+1}$ is a ``weak conditioning'', and to bound certain quantities, this conditioning may be suppressed, leading to easy computations. The archetype of this remark is the following property
\begin{lem}\label{dec} Assume $({\rm H})$. 
There exists $c>0$ such that for $n$ large enough, for any set $A_n$ on ${\cal H}_n$ depending only on the first half part of the trajectories, (that is $\sigma(H_0,H_1, \dots, H_{\floor{n/2}})-$measurable), 
\begin{equation}\label{br}
\QQ(A_n):=`Q_{p_n}^n(A_n \,| \,{\cal H}_n^{2\k{n}+1} )\leq c \, \, `Q_{p_n}^n(A_n).
\end{equation}
\end{lem}
\proof The idea is taken from the proof of Lemma 1 in \cite{JM}~:
\be
\QQ(A_n)&=&\sum_j \QQ(A_n, H_{\floor{n/2}}=j)\\
                                     &=& \sum_j `Q_{p_n}^n(A_n \, |\, H_{\floor{n/2}}=j, H_{n+1}=2\k{n}+1)\,\,  \QQ(H_{\floor{n/2}}=j)\\
&=& \sum_j `Q_{p_n}^n(A_n \, |\, H_{\floor{n/2}}=j)\,\,  \QQ(H_{\floor{n/2}}=j).
\ee
The latter  equality comes from the Markov property of $\H$ under $`Q_{p_n}^n$ that implies that $`Q_{p_n}^n(H_{n+1}=2\k{n}+1|A_n,H_{\floor{n/2}}=j)=`Q_{p_n}^n(H_{n+1}=2\k{n}+1|H_{\floor{n/2}}=j)$. 
It suffices to establish that there exists $c\geq 0$ such that for $n$ large enough, for any $j$
\[ \QQ(H_{\floor{n/2}}=j)\leq c\,\,`Q_{p_n}^n(H_{\floor{n/2}}=j).\]
Write
\[\QQ(H_{\floor{n/2}}=j)=`Q_{p_n}^n(H_{\floor{n/2}}=j)\frac{`Q^n_{p_n}(H_{\floor{n/2}}=2\k{n}+1-j)}{`Q_{p_n}^n(H_{n+1}=2\k{n}+1)};\]
using Lemma \ref{lll}, the last quotient is bounded, uniformly on $j$ and $n\geq 1$.~$\Box$

A simple consequence of Lemma \ref{dec} is the following~: let $X_n$ be a positive random variable that depends only on the first half part of the trajectories, then 
the expectation of $X_n$ under $\QQ$ is bounded by the expectation of $c\,X_n$ under $`Q_{p_n}^n$.

\subsection{Tightness under $`P_n^{{\w},(\k{n})}$}
\label{pnw}
Assume that $(\k{n})$ satisfies $(H)$, that $\S\in \W_n$, and let $\H=\Phi_n(\S$).
Set 
\begin{equation}
h_n(t)=\frac{H_{(n+1)t}}{2\gamma_n}-\frac{t \k{n}}{\gamma_n} \textrm{ for }s,t\in[0,1],
\end{equation}
where $\H$ is assumed to be interpolated between integer points. 
Thanks to formulas \eref{beb}, \eref{sl}, \eref{are},
\begin{equation}\label{concc}
\l|\lambda_n(t)-\lambda_n(s)\r| \leq |h_n(t)-h_n(s)|+2/{\gamma_n}.
\end{equation}
Hence, the tightness of $(h_n)$ under $\QQ$ implies the tightness of $(\lambda_n)$ under $`P_n^{\w,(\k{n})}$.

By symmetry of the random walk under these distributions, we may prove the tightness only on $[0,1/2]$.  By Lemma \ref{dec}, the tightness of $(h_n)$ on $[1/2]$ under $`Q_{p_n}^n$ implies the tightness of $(h_n)$ on $[0,1/2]$ under $\QQ$, and then that of $(\lambda_n)$ on $[0,1/2]$ under $`P_n^{\w,(\k{n})}$. Hence, it suffices to prove the tightness of $(h_n)$ on $[0,1]$ under $`Q_{p_n}^n$ to prove that of $(\lambda_n)$ on $[0,1]$ under $`P_n^{\w,(\k{n})}$.

\begin{comm}
The conditioning by the number of peaks is a strong conditioning on ${\W_n}$. Indeed, ${\cal W}_n^{(\k{n})}$ may have a very small (even exponentially small) probability under $`P_n^\w$ when $\k{n}$ is far away from $n/4$~: no tight bound can be derived using comparison between $`P_n^\w$ and $`P_n^{\w,(\k{n})}$ by just removing the conditioning by $\Lambda_n=\k{n}$. The passage by $`Q_{p_n}^n$ allows to remove this conditioning.
\end{comm}

\subsubsection*{Tightness of the sequence $(\lambda_n)$ under $`P_n^{{\w},(\k{n})}$}
At first sight,  under $`Q_{p_n}^n$, $(h_n)$ is a random walk with the right normalization, and it should converges to the Brownian motion (and then the tightness should follow). However, we were unable to find a reference for this result under the present setting. We will then prove it. \par
In the sub-case where there exists $\delta>0$ such that, for $n$ large enough, $p_n$ satisfies 
$n^{\delta-1}\leq p_n \leq 1-n^{\delta-1}$ then under $`Q_{p_n}^n$, $h_n\dd \w$ in $C[0,1]$~: it is consequence of  Rackauskas, \&  Suquet \cite[Theorem 2]{RS}.  
In this case the tightness holds in a space of H\"older functions, with exponent smaller than $1/2$· When $p_n=o(n^{\delta-1})$ or $1-p_n=o(n^{\delta-1})$, for any $\delta$, $(h_n)$ is not tight in any H\"older space; this may be checked in considering a single normalized step.\medskip

Let $`e>0$ and $\eta>0$ be fixed, and let us prove that for any $n$ large  enough, $`Q_{p_n}^n( \omega_{\delta}(h_n)\geq `e)\leq \eta$ for $\delta$ sufficiently small.
So take a parameter $\delta\in(0,1)$. We have
\begin{equation}\label{edz}
\omega_{\delta}(h_n)\leq 2 \max_{0\le j\le \floor{1/(2\delta)}}  \l(\max_{ I_{j}^{\delta}(n)}h_n-\min_{ I_j^{\delta}(n)}  h_n\r)
\end{equation}
where 
\[I_j^{\delta}(n)=\l[\frac{2j\floor{\delta (n+1)}}{n+1}\wedge 1,\frac{2(j+1)\floor{\delta (n+1)}}{n+1}\wedge 1\r]\] (notice that the length of $I_j^{\delta}(n)$ is larger than $\delta$ for $n$ large enough, and smaller than 3$\delta$). The factor 2 in \eref{edz} simply comes from the splitting up $[0,1]$ into parts. Since, the extremities of the $I_{j}^{\delta}(n)$'s coincide with the discretization points, by the Markov property of  $h_n$,
\begin{equation}\label{inter}
`Q_{p_n}^{n}(\omega_{\delta}(h_n)\geq ` e)\leq (1/(2\delta)+1)\, `Q_{p_n}^{n}\Big(\sup_{ I_1^{\delta}(n)} \l|h_n\r|\geq `e/2\Big).\end{equation}
We need to control the supremun of a random walk, and we then use Lemma \ref{berbound}.  
\begin{lem}\label{zef}
Let $B(p_n)$ be a Bernoulli random variable with parameter $p_n$. There exists $K>0$, such that for any $n\geq 1$,  any $|t|\leq \gamma_n$,
\begin{equation}\label{grs}
`E(e^{t (B(p_n)-p_n)/\gamma_n})\leq \exp\l(\frac{2K}{n} \frac{t^2}{2}\r).
\end{equation}
\end{lem}
\proof There exists $K>0$ such that, for any $|x|\leq 1$, $e^x\leq 1+x+Kx^2$. 
Hence, for any $|t|\leq 1,$ $`E(e^{t (B(p_n)-p_n)})=p_ne^{t(1-p_n)}+(1-p_n)e^{-tp_n}\leq 1+Kt^2(p_n(1-p_n))\leq e^{2Kt^2p_n(1-p_n)}$. Hence, for any $|t|\leq \gamma_n$, \eref{grs} holds  (recall that $(1-p_n)p_n/\gamma_n^2=1/n$).~$\Box$\medskip

Let us end the proof of tightness of $(h_n)$. Since, for $N\in\cro{0,n}$, $h_{N/(n+1)}$ is a sum of $N$ i.i.d. r.v. with the law of $(B(p_n)-p_n)/\gamma_n$, 
\be
`Q_{p_n}^{n}\Big(\sup_{j\leq N} \l|h_{j/(n+1)}\r|\geq `e/2\Big)&\leq& 2\exp\l(-(n+1)`e^2/(16NK)\r) \textrm{ for any }`e\leq 2N \gamma_nK/(n+1). 
\ee
Hence, for $`e>0$, $\delta>0$ fixed, for $n$ large enough,
\be
`Q_{p_n}^{n}(\omega_{\delta}(h_n)\geq ` e)&\leq& (1/(2\delta)+1)\,
 Q_{p_n}^{n}\l( \sup_{0\leq t \leq 3\delta} \l|h_{\floor{(n+1)t}}\r| \geq `e/2 \r)\\
&\leq& 2(1/(2\delta)+1)\,\exp\l(-\frac{`e^2}{16K (3\delta)}\r)
\ee
and this is smaller than any $\eta$ for $\delta$ small enough, and $n$ large enough.~$\Box$ 
\subsubsection*{Tightness of $(s_n)$ under $`P_n^{\w,(\k{n})}$}

Let $\delta\in(0,1)$. First, it suffices to prove the tightness of $(s_n)$ under $`P_n^{\w,(\k{n})}$ on $[1/2]$. Consider formula \eref{gege}. Denote by $G_{i+1}(\H)={\cal T}_{i+1}(\H)- {\cal T}_{i}(\H)$. 
For any $s$ and $t$ such that $ns$ and $nt$ are integer, and $s<t$, we have
\ben\label{grr}
s_n(t)-s_n(s)={\beta_n^{-1}}\sum_{i=H_{ns}}^{H_{nt}-1} (-1)^{i+1} G_{i+1}(\H)+ (-1)^{H_{ns}}\l((ns-{\cal T}_{H_{ns}}(\H))-(nt-{\cal T}_{H_{nt}}(\H))\r).
\een
The range of $s_n$ in a subinterval $I\subset [0,1/2]$ is then a function of the values of $h_n$ in the same interval. Denote by $y_n(s,t)$ the right hand side of \eref{grr}. 
We may control the range of $s_n$ under $`P_n^{\w,(\k{n})}$ by making some computations on $y_n(s,t)$ under $\QQ$, and then, by Lemma \ref{dec}, we may work with $y_n(s,t)$ under $`Q_{p_n}^{n}$.  Under $`Q_{p_n}^{n}$, the variables $G_{i}(\H)$ are geometrical ${\cal G}(p_n)$ distributed ($`P(G_1=k)=p_n(1-p_n)^{k-1}\ind_{k\geq 1}$), and then the computations are simplified.
By Formula \eref{grr}, we have
\be
|y_n(s,t)|&\leq& \l|\sum_{i=H_{ns}}^{H_{nt}-1} {\beta_n^{-1}}(G_{i+2}(\H)-G_{i+1}(\H))\ind_{i\equiv H_{ns} \mod 2}\r|+3{\beta_n^{-1}}\max_{j\leq H_n} G_{j}(\H) 
\ee
where, in the sum, we have packed the variables $G_i(\H)$ per 2. \par
Denote by $\tilde{y}_n(s,t)$ the sum. 
Using that $\max_{k\leq m} G_k$ is a non decreasing function of $m$, we have
\ben
`Q_{p_n}^{n}\l(\frac{\sup_{j\leq H_n} G_{j}(\H)}{\beta_n}\geq `e\r)&\leq& `Q_{p_n}^{n}\l(H_n>4\k{n}\r)+`Q_{p_n}^{n}\l(\frac{\sup_{j\leq 4\k{n}} G_{j}(\H)}{\beta_n}\geq `e\r).
\een
Since $H_n$ is a binomial random variable with parameter $n$ and $p_n=2\k{n}/n$, by the Bienaymé-Tchebichev's inequality, the first term in the right hand side is $O(\k{n}^{-1})$. For the second term, write $`Q_{p_n}^{n}\l(\sup_{j\leq 4\k{n}} G_{j}(\H)\geq `e{\beta_n}\r)=1-(1-`P(G_1\geq`e\beta_n))^{4\k{n}}$. Since $`P(G_1\geq`e\beta_n)=(1-p_n)^{\ceil {`e\beta_n}-1}$, we find that the second term goes to 0 when $n\to +\infty$. \par
It remains to control the variables $\tilde{y}_n(s,t)$. 
Using the Markov property of the random walk $z_k=\sum_{i=1}^{k} {\beta_n^{-1}}(G_{2i}(\H)-G_{2i-1}(\H))$, we get that under $`Q_{p_n}^{n}$, 
\[\sup_{s,t\in I_j^{\delta}(n)}|\tilde{y}_n(s,t)|\sur{=}{(d)} \sup_{s,t\in I_1^{\delta}(n)}|\tilde{y}_n(s,t)|\leq 2\sup_{t\in[0,3\delta]}|\tilde{y}_n(0,t)|.\]
Writing $\widetilde{G}_i$ instead of $G_{2i}-G_{2i-1}$, we have
\[\l\{\sup_{t\leq 3\delta}|\tilde{y}_n(0,t)|\geq \frac{`e}{2}\r\}\subset\l(\l\{H_{3n\delta}>12\k{n}\delta\r\}\bigcup\l\{\sup_{j\leq 6\k{n}\delta}\l|{\beta_n}^{-1}\sum_{i=1}^{j} {\widetilde{G}_{i}(\H)}\r|\geq `e\r\}\r).\]
Once again, by Bienaymé-Tchebichev, $`Q_{p_n}^{n}(\l\{H_{3n\delta}>12\k{n}\delta\r\})=O(1/(\delta \k{n}))$. For the second set in the union, we have to control the maximum of a random walk with increments the variables $\widetilde{G}$. 
\begin{lem}\label{gege2}
Let $G_1$ and $G_2$ be two independent geometrical random variables with parameter $p_n$.
There exists $c>0,c'>0$, such that for any $|t|\leq c'\gamma_n$,
\ben \label{qzd}
`E\l(\exp\l(t\frac{G_1-G_2}{\beta_n}\r)\r)\leq \exp\l(\frac{ct^2}{np_n}\r)
\een
\end{lem}
\proof Write $`E\l(\exp\l(u{(G_1-G_2)}\r)\r)= \frac{p_n^2}{(e^{-u}+p_n-1)(e^u+p_n-1)}.$
The denominator $D$ is equal to $1+(1-p_n)^2-(1-p_n)(e^u+e^{-u})$. There exists a constant $c>0$, such that for any $|u|\leq 1$,  $(e^u+e^{-u})\leq 2+ c u^2$. And then $D\geq p_n^2-c(1-p_n)u^2\geq p_n^2(1-\frac{c(1-p_n)u^2}{p_n^2})\geq p^2_n \exp(-\frac{c}{2}\frac{(1-p_n)u^2}{p_n^2})$ this last equality holds when $\frac{c}{2}\frac{(1-p_n)u^2}{p_n^2}\leq 1$. Hence $`E\l(\exp\l(u{(G_1-G_2)}\r)\r)\leq \exp(\frac{c}{2}\frac{(1-p_n)u^2}{p_n^2})$ for $|u|\leq \sqrt{\frac{2}{c(1-p_n)}}p_n$.
Hence, \eref{qzd} holds for $|t|\leq c'\beta_n p_n/\sqrt{1-p_n}$, and then for $|t|\leq c'\gamma_n$.~$\Box$ \medskip

We end now the proof of tightness for the family $(s_n)$. According to Lemmas \ref{berbound} and \ref{gege2}, for $`e>0,\delta>0$ fixed, for a constant $c''>0$ and $n$ large enough 
\[`Q_{p_n}^{n}\Big(\beta_n^{-1}\sup_{j\leq 6\delta \k{n}} \widetilde{G}_{j}(\H)\geq `e/2\Big)\leq \exp(-`e^2/(c''\delta)).\]
This allows to conclude as for $(h_n)$.~$\Box$

\subsection{Tightness under $`P_n^{{\b},(\k{n})}$} \label{tub}
In this section, $n$ is an even number.
Since $`P_n^{{\b},(k)}$ is the uniform distribution on $\B_n^{(k)}$, it coincides with the conditional law $`P_n^{{\w},(k)}(~.~| S(n)=0)$. 
We first establish a lemma that allows to control the probability of a set under $`P_n^{{\b},(\k{n})}$, by the probability of the same set under $`P_n^{{\w},(k)}$.
\begin{lem}\label{dsad}Assume $(\rm{H})$. There exists a constant $c>0$ such that for $n$ large enough, any $A_n\subset\W_n$,
\[`P_n^{{\b},(\k{n})}(A_n)\leq c \,\beta_n\, `P_n^{\,{\w},(\k{n})}(A_n).\]
\end{lem}
\proof Write
\[`P_n^{{\b},(\k{n})}(A_n)=\frac{`P_n^{{\w},(\k{n})}(A_n, S(n)=0)}{`P_n^{{\w},(\k{n})}(S(n)=0)}\leq \frac{`P_n^{{\w},(\k{n})}(A_n)}{`P_n^{{\w},(\k{n})}(S(n)=0)}.\] 
Now $`P_n^{{\w},(\k{n})}(S(n)=0)=\#\B_n^{(\k{n})} /\#\W_n^{(\k{n})}\sim \sqrt{2/\pi}\beta_n^{-1}$, by the local limit theorem. ~$\Box$\medskip

Since $\gamma_n\to +\infty$, this lemma is interesting only for sets with probability $o(\gamma_n^{-1})$, e.g.~:
\begin{lem}\label{se}
Assume $(\rm{H})$. Let $(\alpha_n)$ be a sequence such that $\alpha_n\to +\infty$.
There exist $c>0$, $c'>0$, such that for $n$ large enough 
\[`P_n^{{\w},(\k{n})}\l(\l|\Lambda_{n/2}-\k{n}/2\r|\geq \alpha_n{\gamma_n}\r)\leq \exp(-c'\alpha_n^2),\]
and then
\[`P_n^{{\b},(\k{n})}\l(\l|\Lambda_{n/2}-\k{n}/2\r|\geq \alpha_n{\gamma_n}\r) \leq c \beta_n\exp(-c'\alpha_n^2).\]
\end{lem}
\proof Assume again that $\H=\Phi_n(\S)$. The variable
$\Lambda_{n/2-1}(\S)=\floor{H_{n/2}/2}$ depends only on the first half of the trajectories. By Lemma \ref{dec}, for any set $I_n$,
\[`P_n^{{\w},(\k{n})}\l(\Lambda_{n/2-1}(\S)\in I_n\r)=\QQ\l(\floor{H_{n/2}/2}\in I_n\r)\leq c'' \, \, `Q_{p_n}^n\l(\floor{H_{n/2}/2}\in I_n\r).\]
Under $`Q_{p_n}^n$, $H_{n/2}$ is a binomial random variable with parameters $n/2$ and $p_n$. Now, using Lemmas \ref{berbound} and \ref{zef} with $I_n=\complement[\k{n}/2-\alpha_n\gamma_n,\k{n}/2+\alpha_n\gamma_n]$, we get the first assertion. The second assertion is a consequence of Lemma \ref{dsad}.~$\Box$\medskip

Consider $\bar{\W_n}$ the set of simple walks satisfying
\[\bar{\W_n}=\{\S: \S\in {\W}_{n}, \l|\Lambda_{n/2}-\k{n}/2\r|\leq \gamma_n^{5/4} \}.\]
Since by Lemma \ref{se}, $`P_n^{{\b},(\k{n})}(\bar{\W_n})\to 1$, we will from now on concentrate on these trajectories. We stress on the fact that $\gamma_n^{5/4}=o(\k{n})$. Assume that the following lemma is proved. 
\begin{lem}\label{condb} Assume $({\rm H})$. 
There exists a constant $c>0$ such that for $n$ large enough, for any subset $A_n\subset \bar{\W_n}$ depending only on the first half of the trajectories, 
\begin{equation}\label{zr}
`P_n^{{\b},(\k{n})}(A)= `P_n^{{\w},(\k{n})}(A_n\,| \,S(n)=0 )\leq c \, \, `P_n^{{\w},(\k{n})}(A_n).
\end{equation}
\end{lem}
This lemma, very similar to Lemma \ref{dec}, allows to obtain the tightness of $(s_n,\lambda_n)$ under $`P_n^{{\b},(\k{n})}$ from that under $`P_n^{{\w},(\k{n})}$; proceed as follows. 
By symmetry of bridges under $`P_n^{{\b},(\k{n})}$, it suffices to prove the tightness on $[0,1/2]$. Since $`P_n^{{\b},(\k{n})}(\bar{\W_n})\to 1$, we restrict our study to the trajectories of $\bar{\W_n}$. By \eref{zr}, the tightness of $(s_n,\lambda_n)$ under  $`P_n^{{\w},(\k{n})}$ on $[0,1]$ (and then on $[0,1/2]$) implies that under $`P_n^{{\b},(\k{n})}$ on $[0,1/2]$. It only remains to prove Lemma \ref{condb}.\medskip

\noindent\bf Proof of Lemma \ref{condb} \rm 
First, for any $A_n$, depending only on the first half  of the trajectories, 
\ben\label{aza}
`P_n^{{\w},(\k{n})}(A_n |S(n)=0)&=&\sum_{(l,x,a)}`P_n^{{\w},(\k{n})}(A_n\,|\,(\Lambda_{n/2},S({n/2}),\Delta S_{n/2})=(l,x,a))\\&&\times\,`P_n^{{\w},(\k{n})}((\Lambda_{n/2},S({n/2}),\Delta S_{n/2})=(l,x,a)~|~S(n)=0)
\een
where the summation is taken on all possible triples $(l,x,a)$. Indeed, under $`P_n^{\w}$, the sequence $Y_k=(\Lambda_k,S(k),\Delta S_k)$ is a Markov chain, and also under $`P_n^{\w,(\k{n})}$. Then write $`P_n^{{\w},(\k{n})}(A_n |S(n)=0)= \sum_{(l,x,a)} `P_n^{{\w},(\k{n})}(S(n)=0 | A_n, Y_{n/2}=(l,x,a))`P_n^{{\w},(\k{n})}(Y_{n/2}=(l,x,a))/`P_n^{{\w},(\k{n})}(S(n)=0)$, then the first conditioning by $A_n$ may be deleted, by Markov. A trite computation leads to the result. \par
Now, assume that $A_n\subset \bar{\W_n}$. The summation in \eref{aza} can be done on the triples $(l,x,a)$, such that $l\in J_n:= \cro{\k{n}/2-\gamma_n^{5/4},\k{n}/2+\gamma_n^{5/4}}$, $x\in \cro{-n/2,n/2}$, $a\in\times\{\u,\d\}$.
To end the proof, we check that there exists $c>0$, valid for any $(l,x,a)\in J_n\times\cro{-n/2,n/2}\times\{\u,\d\}$, and $n$ large enough, such that 
\ben \label{zzer}
`P_n^{{\w},(\k{n})}\l(Y_{n/2}=(l,x,a)~|~S(n)=0\r)\leq c\,`P_n^{{\w},(\k{n})}\l(Y_{n/2}=(l,x,a)\r).
\een
We choose to condition by the last increment of $\S$ in $\cro{0,n/2}$ for computation reasons.\par
For any $a\in\{\u,\d\}$, denote $T_{-a}=T_{\u a}\cup T_{\d a}$ and similar notation for $T_{a-}$ and $T_{-\,-}$. \par
$\bullet$ Case $a=d$. In this case, the left hand side of \eref{zzer} equals 
\[\frac{\#T_{-\,\d}^{l}(n/2,x,0)\#T_{-\,-}^{\k{n}-l}(n/2,x,0)}{\# \B_n^{(\k{n})}}\] and the right hand side
\[\frac{\#T_{-\,\d}^{l}(n/2,x,0)\#\W_{n/2}^{(\k{n}-l)}}{\#\W_n^{(\k{n})}}.\]
Since for any $c>0$, the application $x\mapsto\Gamma(x)/\Gamma(x-c)$ is $\log$-concave, 
\[\#T_{-\,-}^{\k{n}-l}(n/2,x,0)=\binom{\frac{n/2-x}{2}}{\k{n}-l}\binom{\frac{n/2+x}2}{\k{n}-l}\leq \binom{\floor{n/4}}{\k{n}-l}\binom{\ceil{n/4}}{\k{n}-l}:=g_{n,l}.\]
Using that  
\[\frac{\#\B_n^{(\k{n})}}{\#\W_n^{(\k{n})}}=\frac{p_n \#\B_n^{(\k{n})}p_n^{2\k{n}}(1-p_n)^{n-2\k{n}}}{\#\W_n^{(\k{n})}p_n^{2\k{n}+1}(1-p_n)^{n-2\k{n}}}\sim p_n\frac{\sqrt{2\pi}\gamma_n}{2\pi \gamma_n^2/2}\]
To prove \eref{zzer} when $a=d$ it suffices to prove that 
\[\limsup_n \max_{l\in J_n}\frac{\gamma_n}{p_n}\frac{g_{n,l}}{\#\W_{n/2}^{(\k{n}-l)}}<+\infty\]
(notice that $\floor{n/4}+\ceil{n/4}=n/2$ for any even $n$.)
We have
\[\frac{\gamma_n}{p_n}\frac{g_{n,l}}{\#\W_{n/2}^{\k{n}-l}}=\gamma_n\frac{\binom{\floor{n/4}}{\k{n}-l}\binom{\ceil{n/4}}{\k{n}-l}}{\binom{n/2}{2(\k{n}-l)}\l(\frac{n/2+1}{2\k{n}-2l+1}p_n\r)}\]
The (second) parenthesis in the denominator is bounded. It remains to prove that $\limsup_n \max_{l\in J_n} G_n(l)$ is bounded, for 
\[G_n(l):= \frac{\gamma_n\binom{\floor{n/4}}{\k{n}-l}\binom{\ceil{n/4}}{\k{n}-l}}{\binom{n/2}{2(\k{n}-l)}}.\]
We will prove this assertion by showing that for any sequence $(l_n)$ of integers, that satisfies $l_n\in J_n$ for any $n$,   $(G_n(l_n))$ converges to a constant that does not depend on $(l_n)$. This allows to conclude, since one may take the sequence $(l_n)$ s.t. $G_n(l_n)$ maximizes $G_n(l)$ on $J_n$ for any $n$. Set $\rho_{n}=4(\k{n}-l_n)/n$. Since $l_n\in J_n$,  for any $n$, $\rho_n\in(0,1)$.Now, one checks easily that 
\[G_n(l_n)=\frac{\gamma_n`P\l(B(\floor{n/4},\rho_{n})=\k{n}-l_n\r)`P\l(B(\ceil{n/4},\rho_{n})=\k{n}-l_n\r)}{`P\l(B({n/2},\rho_{n})=2(\k{n}-l_n)\r)}. \] 
Since $\k{n}-l_n\sim \k{n}/2$, by the central local limit theorem, it converge to $2/\sqrt{\pi}$.\\
$\bullet$ Case $a=\u$. In this case, the left hand side of \eref{zzer} equals 
\[\frac{\#T_{-\,\u}^{l}(n/2,x,0)(\#T_{u\,-}^{\k{n}-l}(n/2,x,0)+\#T_{d\,-}^{\k{n}-l-1}(n/2,x,0))}{\# \B_n^{(\k{n})}}\] and the right hand side
\[\frac{\#T_{-\,\d}^{l}(n/2,x,0)(\#\W_{n/2,\u}^{(\k{n}-l)}+\#\W_{n/2,\d}^{(\k{n}-l-1}))}{\#\W_n^{(\k{n})}},\]
where $\W_{n,a}^k $ is the set of trajectories $\S$ with $k$ peaks with $\Delta S_1=a$.\par
Once again, it suffices to check that the quotient
\[\#T_{u\,-}^{\k{n}-l}(n/2,x,0)+\#T_{d\,-}^{\k{n}-l-1}(n/2,x,0)=\binom{1+\frac{n/2-x}{2}}{k-l}\binom{\frac{n/2+x}2-1}{\k{n}-l-1}\]divided by 
\[\#\W_{n/2,\u}^{\k{n}-l}+\#\W_{n/2,\d}^{(\k{n}-l-1)}=\binom{n/2+1}{2(\k{n}-l)}\]
is bounded by $c \, {\# \B_n^{(\k{n})}}/{\#\W_n^{(\k{n})}}$. The same arguments leads to the same conclusion.$\Box$

\subsection{Tightness under $`P_n^{\,{\m},(\k{n})}$}\label{bml}
\subsubsection*{The case $n$ even}
\begin{figure}[htbp]
\psfrag{a}{An excursion-type path}\psfrag{B}{A bridge $B$}\psfrag{P}{$\Psi_{2N}(B)$}
\psfrag{0}{0}\psfrag{t_1}{$t_1$}\psfrag{t_2}{$t_2$}\psfrag{t_3}{}\psfrag{t_4}{}
\centerline{\includegraphics[height=5cm,width=13cm]{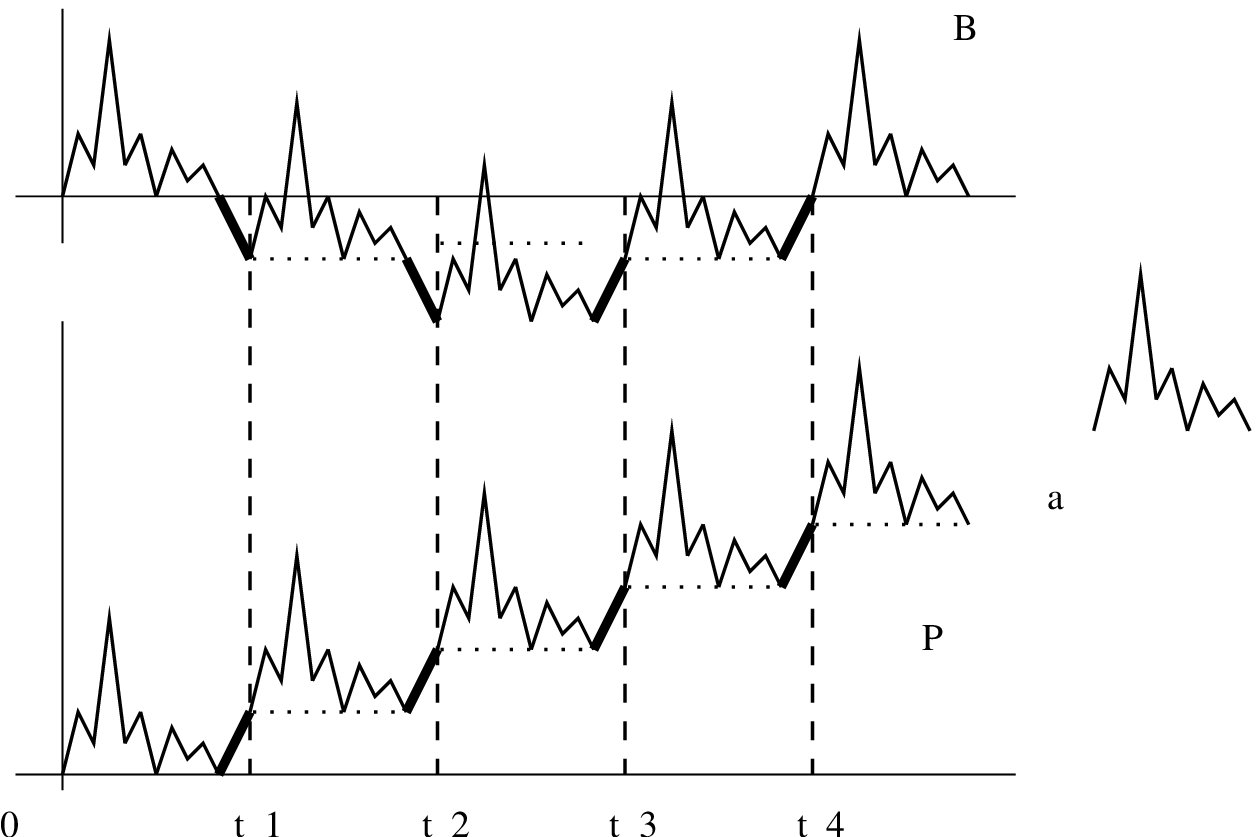}}
\captionn{\label{Psi} Synthetic description of $\Psi_{2N}$. The application $\Psi_{2N}$ turns over each increment corresponding to a reaching time of a negative position. 
The application $\Psi_{2N}^{-1}$ turns over the last increments reaching each position $x\in \cro{1,Z_{2N}/2}$ ($Z_{2N}$ is even).}
\end{figure}

Assume first that $n=2N$ is even.
We recall a bijection $\Psi_{2N}:\B_{2N}\to \M_{2N}$, illustrated on Figure \ref{Psi}, that maps $\B_{2N}^{(k)}$ on $\M_{2N}^{(k)}$, and that moreover preserves sufficiently the trajectories, to prove that the tightness of $(s_{2N},\lambda_{2N})$  under $`P_{2N}^{\,{\b},(\k{2N})}$ yields that under $`P_{2N}^{\,{\m},(\k{2N})}$.\medskip 

The application $\Psi_{2N}:\B_{2N}\mapsto\W_{2N}$ (we will see later that $\Psi_{2N}(\B_{2N})=\M_{2N}$) is defined as follows.
Let $\S\in \B_{2N}$ and $m=\min \S\leq 0$ its minimum. 
For $j\in \{1,\dots,-m\}$, let $t_j=\tau_{-j}(\S)$ the reaching time of $-i$ by $\S$. Write $I_{\S}=\{t_j,j\geq 1\}$. Notice that when $m=0$,  $I_{\S}=\varnothing$. 

The trajectory $\Psi_{2N}(\S)=\Z=(Z_i)_{i=0,\dots,n}$ is defined by $Z_0=0$ and ~:
\[\Delta Z_i=\left\{
\begin{array}{ll}
\Delta S_i &\textrm{ if } i\notin I_{\S},\\
-\Delta S_i=+1 &\textrm{ if } i\in I_{\S}.
\end{array}
\right.
\]
\begin{pro}\label{psin}
For any even ${2N}$, $\Psi_{2N}$ is a bijection from $\B_{2N}$ onto $\W_{2N}$; moreover, for any $k$, its restriction to $\B_{2N}^{(k)}$ is a bijection onto $\M_{2N}^{(k)}$ that preserves the peak positions. 
\end{pro}

\proof First, it is easy to see that if $\Z=\Psi_{2N}(\S)$, for any $i\leq {2N}$, 
\begin{equation}\label{aar}
Z_i=S(i)+2\min_{j\leq i}S(j).
\end{equation}
Hence, $\Psi_{2N}(\B_{2N})\subset \M_{2N}$. Since $\Psi_{2N}$ is clearly an injection, the first assertion of the Proposition follows $\#\B_{2N}=\#\M_{2N}$. Since $\Psi_{2N}$ does not create or destroys any peaks, or even changes the position of the peaks, the restriction of $\Psi_{2N}$ onto $\B_{2N}^{(k)}$ is a bijection onto $\Psi_{2N}(\B_{2N}^{(k)})\subset \M_{2N}^{(k)}$. The equality $\#\B_{2N}^{(k)}=\#\M_{2N}^{(k)}$ suffices then to conclude. ~$\Box$

By Proposition \ref{psin}, the tightness of $(\lambda_{2N})$ under $`P_{2N}^{\,{\m},(\k{2N})}$  is  a  consequence of that under $`P_{2N}^{\,{\b},(\k{2N})}$. For $(s_{2N})$, \eref{aar} implies that the modulus of continuity of the non normalized trajectories are related by 
\[\omega_{\Z}(\delta)\leq 3\omega_{\S}(\delta)\textrm{ for any }\delta\in[0,n]\]and then, the tightness of $(s_{2N})$ under $`P_{2N}^{\,{\m},(\k{2N})}$  follows that under $`P_{2N}^{\,{\b},(\k{2N})}$.

\subsubsection*{The case $n$ odd}
The case $n=2N+1$ odd is very similar.
There exists a bijection $\Psi_{2N+1}$ between  $\widetilde{\B}_{2N+1}$ and $\M_{2N+1}$ where $\widetilde{\B}_{2N+1}$ is the subset of $\W_{2N+1}$ of trajectories ending at position $+1$.\begin{figure}[htbp]
\psfrag{a}{An excursion-type path}\psfrag{B}{A bridge $\tilde{B}$}\psfrag{P}{$\Psi_{2N+1}(\tilde{B})$}
\psfrag{0}{0}\psfrag{t_1}{$t_1$}\psfrag{t_2}{$t_2$}\psfrag{t_3}{}\psfrag{t_4}{}
\centerline{\includegraphics[height=5cm,width=13cm]{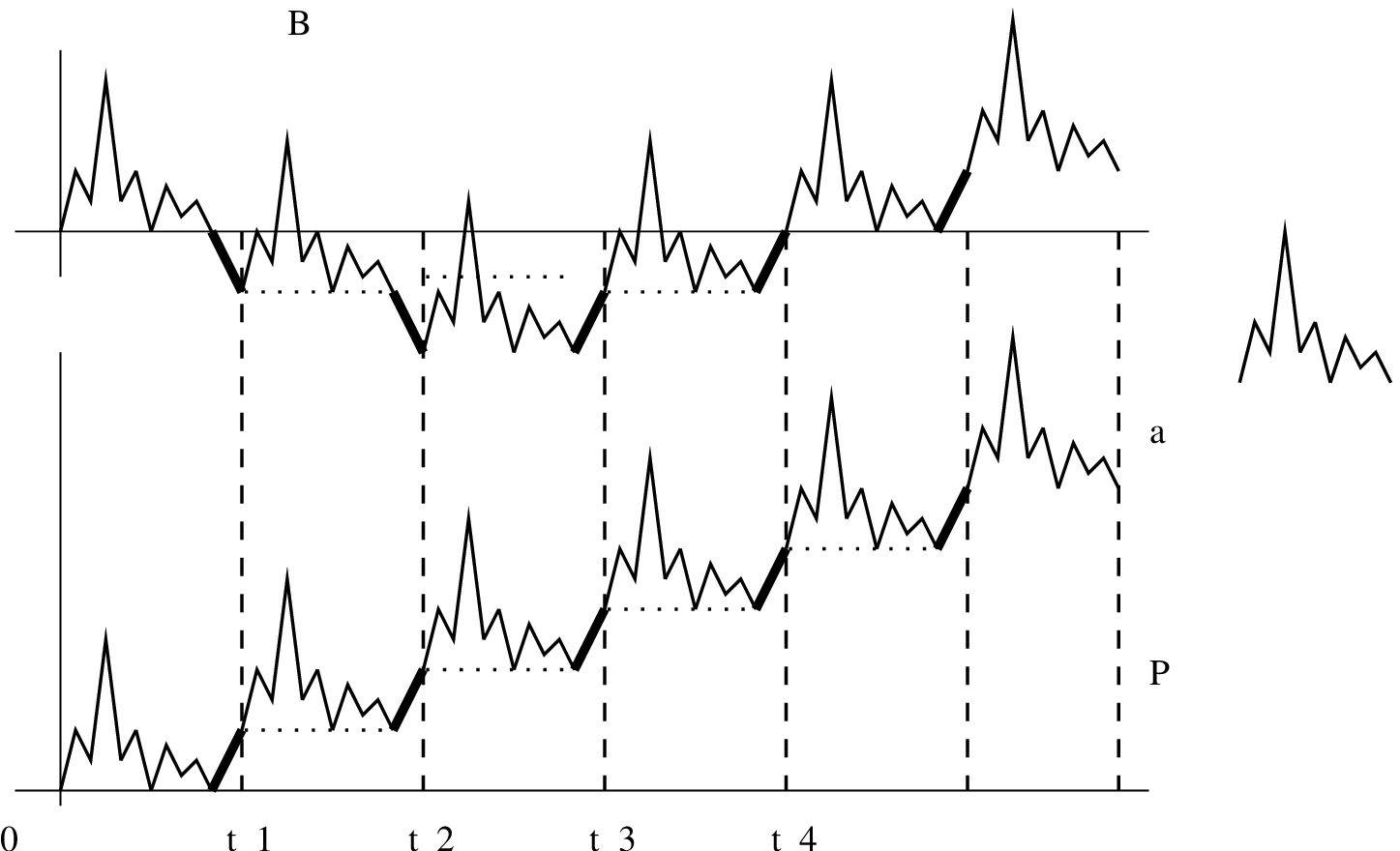}}
\captionn{\label{Psi2} Synthetic description of $\Psi_{2N+1}$. The application $\Psi_{2N+1}$ turns over each increment corresponding to a reaching time of a negative position. The application $\Psi_{2N+1}^{-1}$ turns over the last increments reaching each position $x\in \cro{1,\floor{Z_{2N+1}/2}}$ ($Z_{2N+1}$ is odd).}
\end{figure}
The application $\Psi_{2N+1}$ has the same properties as $\Psi_{2N}$ to conserve the peak positions, and  the set $\tilde{B}_{2N+1}^{(\k{2N+1})}:=\widetilde{\B}_{2N+1}\cap\W_{2N+1}^{(\k{2N+1})}$ is sent on $\M_{2N+1}^{(\k{2N+1})}$. To conclude, we need a tightness result for the uniform distribution on $\tilde{B}_{2N+1}^{(\k{2N+1})}$. But the result of Section \ref{tub} regarding $\B_{2N}^{(\k{2N})}$ maybe generalized to $\B_{2N+1}^{(\k{2N+1})}$. \par
Hence, $(s_{2N+1},\lambda_{2N+1})$ is tight under $`P_{2N+1}^{\,{\m},(\k{2N+1})}$, and then we may conclude that  $(s_{n},\lambda_{n})$ is tight under $`P_{n}^{\,{\m},(\k{n})}$.

\subsection{Tightness under $`P_n^{\,{\e},(\k{n})}$}
Here $n=2N$ is an even number.
Consider 
\be
\check{B}_{2N+1}&=&\{\S, \S\in\W_{2N+1}, S({2N+1})=-1\},\\
\check{E}_{2N+1}&=&\{\S, \S\in\W_{2N+1}, S(j)\geq 0 \textrm{ for any }j\in\cro{0,2N}, S({2N+1})=-1\}
\ee 
and $\check{B}_{2N+1}^{(K)}=\check{B}_{2N+1}\cap \W_{2N+1}^{(K)}$, $\check{E}_{2N+1}^{(K)}=\check{E}_{2N+1}\cap \W_{2N+1}^{(K)}$. Informally, $\check{E}_{2N+1}$  (resp. $\check{E}_{2N+1}^{(K)}$) are Bernoulli excursion from $\E_{2N}$ (resp. with $K$ peaks)  with an additional ending $d$-step, and  $\check{B}_{2N+1}$ and $\check{B}_{2N+1}^{(K)}$ are trajectories ending at $-1$ (resp. with $K$ peaks).\par
Consider the application 
\[\app{\R}{\check{E}_{2N+1}\times\cro{0,2N}}{\check{B}_{2N+1}}{(\S,\theta)}{\R(\S,\theta)=\S^{(\theta)}=(S^{(\theta)})_{i=0,\dots,2N+1}}\]
defined by 
\[\Delta S^{(\theta)}_{k}=\Delta S_{k+\theta \mod 2N+1}\]
or equivalently $S^{(\theta)}_k=S({k+\theta\mod_{2N+1}})-S({\theta})-\ind_{k+\theta>2N+1}$. Informally, $\S\mapsto \S^{(\theta)}$ exchanges the $\theta$ first steps of $\S$ with the last $2N+1-\theta$'s ones. 
The application $\R$ is a bijection between $\check{E}_{2N+1}\times\cro{0,2N}$ and $\check{B}_{2N+1}$: this is the so-called cyclical lemma attributed to Dvoretzky-Motzkin, or Kemperman, or Otter see Pitman \cite[Chapter 5]{PITT} and Bertoin \& al. \cite{BCP}. 
The peaks positions of $\S^{(\theta)}$  are obtained from that of $\S$ by a shift of $-\theta \mod 2N+1$, and $\#\S^{(\theta)}_{\wedge}=\#\S_{\wedge}-1$ iff $\theta\in \S_{\wedge}$ (if $\theta\notin \S_{\wedge}$ then $\#\S^{(\theta)}_{\wedge}=\#\S_{\wedge}$). \\
For any $\S\in \W_{2N+1}$, set $\Xi_0(\S)=0$, and for any $k\leq 2N-\#\S_{\wedge}$, 
\[\Xi_m(\S):=\min\{j, j\geq \Xi_{m-1}(\S), j\leq 2N, j\notin \S_{\wedge}\},\]
the successive non-peak positions of $\S$ in $\cro{0,2N}$. 
Consider the application 
\[\app{\widehat{\R}}{\check{E}_{2N+1}^{(K)}\times\cro{0,2N-K}}{\check{B}_{2N+1}^{(K)}}{(\S, \ell)}{\widehat{\R}(\S,\ell)=\S^{(\Xi_\ell(\S))}}.\] 
\begin{pro}For any $K\leq N$, any $N\geq 0$, the application $\widehat{\R}$ is a bijection from $\check{E}_{2N+1}^{(K)}\times\cro{0,2N-K}$  onto $\check{B}_{2N+1}^{(K)}$.
\end{pro}
\proof It is a consequence of the two following points: for any $\S$, $m\mapsto \Xi_m(\S)$ is a bijection from $\cro{0,\#\cro{0,2N}\setminus \S_{\wedge}}$ onto $\cro{0,2N}\setminus \S_{\wedge}$, and $\R$ is a bijection.~$\Box$~\medskip

Consider $(\S,\theta)$ in $\check{E}_{2N+1}^{(K)}\times\cro{0,2N-K}$; for any $u$,
\ben\label{trze}
\sup_{|m_1-m_2|\leq u}|S({m_1})-S({m_2})|&\leq& 2 \sup_{|m_1-m_2|\leq u}|S^{(\theta)}({m_1})-S^{(\theta)}({m_2})|\\
\nonumber\sup_{|m_1-m_2|\leq u}\l|(\Lambda_{m_1}-\Lambda_{m_2})(\S)-\frac{m_1-m_2}{2N}\k{2N}\r| &\leq& 2 \sup_{|m_1-m_2|\leq u}\l|(\Lambda_{m_1}-\Lambda_{m_2})(\S^{(\theta)})-\frac{m_1-m_2}{2N}\k{2N}\r|.
\een
Endow $\check{E}_{2N+1}^{(K)}\times\cro{0,2N-K}$ with the uniform distribution and consider a random element $(\S,\theta)$ under this law ($\S$ is then uniform on $\check{E}_{2N+1}^{(K)}$). By the last proposition, $\S^{(\theta)}$ is uniform on $\check{B}_{2N+1}^{(K)}$. 
By \eref{trze},  we have
\[`P_{2N}^{e,(\k{2N})}(\omega_{\delta}(s_{2N})\geq `e)\leq `P_{2N+1}^{\w,(\k{{2N}})}(2\omega_{\delta}(s_{2N})\geq `e | S(2N+1)=-1)\]
and the same result holds for $\lambda_{2N}$. Once again  the result of Section \ref{tub} concerning $`P_{2N}^{\b,(\k{2N})}=`P_{2N}^{\w,(\k{2N})}(. | S({2N=0}))$ can be generalized to $`P_{2N}^{\w,(\k{2N})}(. | S(2N+1)=-1)$.  ~$\Box$

\small

\subsubsection*{Acknowledgments} We would like to thank Mireille Bousquet-Mélou who pointed many references, and for helpful discussions.

\end{document}